\documentclass{amsart}
\usepackage[dvips]{graphics} 
\input{diagrams}
\diagramstyle[PostScript=dvips]
\newcommand{\functor}[1]{\mbox{\scshape {#1}}}
\newcommand{\stack}[1]{\mathfrak {#1}}
\newcommand{\scheme}[1]{{\mathcal {#1}}}

\newcommand{\schememg}{\scheme{M}_g}

\newcommand{\mgbar}{\overline{\stack{M}}_g}
\newcommand{\mgn}{\stack{M}_{g,n}}
\newcommand{\schememgn}{\scheme{M}_{g,n}}
\newcommand{\schememgnbar}{\overline{\scheme{M}}_{g,n}}
\newcommand{\mgnbar}{\overline{\stack{M}}_{g,n}}

\newcommand{\sheafhom}{\ch\kern-.15em om}

\newcommand{\sheafext}{\ce\kern-.15em xt}

\newcommand{\sheaftor}{\ct\kern-.15em or}

\newcommand{\spec}{\operatorname{Spec\,}}

\newcommand{\jac}{\operatorname{Jac\,}}

\newcommand{\tensor}{\otimes}
\newcommand{\cross}{\times}
\newcommand{\maxid}{\mathfrak m}
\newcommand{\irightarrow}{\rTo^{\sim}}

\newcommand{\ck}{{\mathcal K}}
\newcommand{\cl}{{\mathcal L}}
\newcommand{\ce}{{\mathcal E}}
\newcommand{\cf}{{\mathcal F}}
\newcommand{\co}{{\mathcal O}}
\newcommand{\cc}{{\mathcal C}}

\newcommand{\ch}{{\mathcal H}}

\newcommand{\cm}{{\mathcal M}}

\newcommand{\cn}{{\mathcal N}}

\newcommand{\ct}{{\mathcal T}}

\newcommand{\fp}{{\mathfrak p}}

\theoremstyle{plain}
\newtheorem{theorem}{Theorem}[subsection]
\newtheorem{proposition}[theorem]{Proposition}
\newtheorem{lemma}[theorem]{Lemma}

\newtheorem{fact}[theorem]{Fact}

\theoremstyle{remark}

\newtheorem{nb}[theorem]{Note}

\newtheorem{example}[theorem]{Example}

\theoremstyle{definition}
\newtheorem{defn}[theorem]{Definition}

\newcommand{\stackzgn}{\stack{C}_{g,n}}

\newcommand{\stackcgn}{\stack{C}_{g,n}}
\newcommand{\stacksgrnmkbar}{\overline{\stack{S}}^{1/r,\bm}_{g,n}(\ck)}

\newcommand{\cgnvbt}{stack{C}_{g,n}}

\newcommand{\stackmgnbar}{\overline{\stack{M}}_{g,n}}
\newcommand{\purerg}{{\functor{Pure}}_{r,g}}

\newcommand{\sym}{\mathbf{Sym}}

\newcommand{\sgrnkbar}{\overline{\stack{S}}_{g,n}^{1/r}(\ck)}

\newcommand{\sgrnmkbar}{\overline{\stack{S}}_{g,n}^{1/r,\mathbf{m}}(\ck)}

\newcommand{\sgr}{{\scheme S}^{1/r}_g}

\newcommand{\sgrnm}{{\scheme S}^{1/r,\mathbf m}_{g,n}}

\newcommand{\sgrbar}{\overline{\scheme S}^{1/r}_{g}}

\newcommand{\sgrnmbar}{\overline{\scheme S}^{1/r,\mathbf m}_{g,n}}

\newcommand{\stacksgr}{{\stack S}_{g}^{1/r}}
\newcommand{\stacksgrbar}{\overline{{\stack S}}_{g}^{1/r}}
\newcommand{\stacksgrnmbar}{\overline{\stack S}^{1/r,\mathbf m}_{g,n}}
\newcommand{\stacksgrnm}{{\stack S}^{1/r,\mathbf m}_{g,n}}

\newcommand{\stackmgbar}{\overline{\stack{M}}_g}

\newcommand{\bm}{\mathbf{m}}          % 

\newcommand{\cv}{c^{virt}}            % The virtual class (Agrees with Tyler)
\begin{document}

\setcounter{tocdepth}{2}

\title[Geometry of spin curve moduli]{Geometry of the moduli of
higher spin curves}

\author{Tyler J. Jarvis}
\address{Department of Mathematics\\ Brigham Young University\\Provo, UT 84602}
\email{jarvis@math.byu.edu}
%\urladdr{http://www.math.byu.edu/~jarvis}
\thanks{This material is based in part on work supported by the National
Science Foundation under Grant No.~DMS-9501617} 
\date{\today}

\begin{abstract}
  
  This article treats various aspects of the geometry of the moduli
  $\sgr$ of $r$-spin curves and its compactification $\sgrbar$.
  Generalized spin curves, or $r$-spin curves, are a natural
  generalization of $2$-spin curves (algebraic curves with a
  theta-characteristic), and have been of interest lately because of
  the similarities between the intersection theory of these moduli
  spaces and that of the moduli of stable maps.  In particular, these
  spaces are the subject of a remarkable conjecture of E. Witten
  relating their intersection theory to the Gelfand-Dikii ($KdV_r$)
  heirarchy.  There is also a $W$-algebra conjecture for these spaces
  \cite{jkv} generalizing the the Virasoro conjecture of quantum
  cohomology.
   
  For any line bundle $\ck$ on the universal curve over the stack of
  stable curves, there is a smooth stack $\stack{S}^{1/r}_{g,n}(\ck)$
  of triples $(X,\cl,b)$ of a smooth curve $X$, a line bundle $\cl$ on
  $X$, and an isomorphism $b:\cl^{\otimes r} \rTo \ck$.  In the
  special case that $\ck = \omega$ is the relative dualizing sheaf,
  then $ \stack{S}^{1/r}_{g,n} (\ck)$ is the stack $\stack{S}^{1/r}_{g,n}$
  of $r$-spin curves.
  
  We construct a smooth compactification
  $\overline{\stack{S}}^{1/r}_{g,n}(\ck)$ of the stack
  $\stack{S}^{1/r}_{g,n}(\ck)$, describe the geometric meaning of its
  points, and prove that it is projective.
  
  We also prove that when $r$ is odd and $g>1$, the compactified stack
  of spin curves $\stacksgrbar$ and its coarse moduli space $\sgrbar$
  are irreducible, and when $r$ is even and $g>1$, $\sgrbar$ is the
  disjoint union of two irreducible components.  We give similar
  results for $n$-pointed spin curves, as required for Witten's
  conjecture, and also generalize to the $n$-pointed case the
  classical fact that when $g=1$, $\overline{\scheme{S}}^{1/r}_1$ is
  the disjoint union of $d(r)$ components, where $d(r)$ is the number
  of positive divisors of $r$. These irreducibility properties are
  important in the study of the Picard group of $\stacksgrbar$
  \cite{jarvis:picard}, and also in the study of the cohomological
  field theory related to Witten's conjecture
  \cite{jkv,witten:r-spin-conj}.

\end{abstract}

\maketitle

\tableofcontents

\section{Introduction}

In this paper we study the geometry of the stack $\stacksgr$ of 
$r$-spin curves (smooth algebraic curves $X$ with a line bundle 
(invertible coherent sheaf) $\cl$ and an isomorphism from the 
$r$th tensor power $\cl^{\tensor r}$ to the cotangent bundle 
$\omega_X$), as well as the coarse moduli $\sgr$, and also the 
compactification $\stacksgrbar$ of $\stacksgr$ by \emph{coherent 
nets of $r$th roots} over the stack of stable curves 
$\stackmgnbar$.

Spin curves  provide  a finite cover of the moduli of curves, 
distinct (in the case of $g>1$) from the moduli of curves with 
level-$r$ structure.  Nor is it the same as the non-Abelian level 
structures that have been studied recently by Looijenga 
\cite{looijenga:level} and Pikaart-de~Jong \cite{pikaart-dejong}.  
In the case of $g=n=1$, however, this space is a disjoint union 
of the modular curves $Y_1(d)$ for all $d$ dividing $r$.  

This moduli space has many similarities to the moduli of  stable 
maps, including the existence of classes analogous to 
Gromov-Witten classes and an associated cohomological field 
theory \cite{jkv}.  But it is not the moduli of stable maps into 
any variety \cite[\S5.1]{jkv}.  It can, however, be described as 
a closed and open substack of the stack of stable maps of 
balanced twisted curves in the sense of Abramovich and Vistoli 
\cite{abramovich-vistoli:complete-moduli} into a Deligne-Mumford 
stack (see \cite{abramovich-jarvis}). 

These moduli spaces are especially interesting because of a 
conjecture of E. Witten, first described in 
\cite{witten:N-matrix-model,witten:r-spin-conj}, relating the 
intersection theory on the (compactified) moduli space of 
$r$-spin curves and Gelfand-Dikii hierarchies of order $r$. This 
conjecture is a generalization of an earlier conjecture of his, 
which was proved by Kontsevich (see 
\cite{kontsevich:thesis,kontsevich:matrix-airy} and 
\cite{looijenga:kontsevich}).  As in the case of Gromov-Witten 
theory, one can construct a potential function from the 
intersection numbers of a cohomology class $\cv$ with the 
tautological $\psi$ classes associated to the universal curve.  
Witten conjectures that the potential of the theory corresponds 
to the tau-function of the order-$r$ Gelfand-Dikii ($KdV_r$) 
hierarchy.   In genus zero, the Witten conjecture is true 
\cite{jkv}, but in higher genus it remains unproven. 

Construction of the moduli stack $\stacksgr$, and its 
compactification $\stacksgrbar$, was done in 
\cite{cornalba:theta} for $r=2$.  For $r 
\geq 2$ several compactifications were constructed in 
\cite{jarvis:spin}.  Unfortunately, none of those 
compactifications is smooth for general $r$, although the space 
$\mbox{\sc Pure}_{r,g}$ of \cite{jarvis:spin} is smooth when $r$ 
is prime.  This article has two main goals:  first, to construct  
a geometrically meaningful compactification which is smooth for 
all $r$, and second, to describe the irreducible and connected 
components of this compactification. 

\subsection{Overview and Outline of the Paper}
\subsubsection{Construction of the Compactification}

For any line bundle $\ck$ on the universal curve $\stack{C}$ over 
the stack $\mgbar$ of stable curves we will denote by 
$\stack{S}^{1/r}_g (\ck)$ the stack of triples $(X,\cl,b)$, where 
$X$ is a smooth curve of genus $g$, $\cl$ is a line bundle on 
$X$, and $b$ is an isomorphism $b:\cl^{\otimes r} \rTo \ck$.  
Thus the stack $\stacksgr$ of smooth spin curves is simply 
$\stacksgr (\omega)$.  We will compactify the stack 
$\stacksgr(\ck)$ for all $r$, $g$, and $\ck$, by considering the 
stack of \emph{coherent nets of $r$th roots}, as given in 
Definition~\ref{coherent}.  Although the definitions of families 
of these nets and the proofs of smoothness and other properties 
are relatively technical, the intuitive ideas are 
straightforward. Moreover, for many applications (for example, 
those in \cite{jarvis:picard,jkv}) one only needs to know the 
definitions for geometric points of the stack, which are much 
simpler, the structure of the universal deformation, as described 
in Theorem~\ref{univ}, and the description of the irreducible 
components, as given in Theorem~\ref{disjoint}. 

The basic idea is simply that, up to isomorphism, an $r$th root 
$\cl$ of a bundle $\ck$ completely determines a whole collection 
(net) of roots $\cl^{\otimes d}$ for each $d$ dividing $r$.  On 
the boundary, the limit of $\cl$ is a rank-one, torsion-free 
sheaf, but the degeneration of $\cl$ is no longer enough to 
completely specify a collection of $d$th roots, if $r$ is not 
prime, so that structure must be made explicit. Moreover, since 
there are generally too many torsion-free sheaves in the limit of 
a degenerating invertible sheaf, we also need to rigidify the 
structure by explicitly specifying the isomorphisms $\cl^{\otimes 
r} 
\irightarrow \omega_X$ and $(\cl^{\otimes d})^{\otimes e} 
\irightarrow \cl^{\otimes de}$ (or comparable homomorphisms 
which are almost isomorphisms over a non-smooth curve). 

Thus the geometric objects of interest will be a collection 
$\{\ce_d, c_{d,d'}\}$, called a net, of a rank-one torsion-free 
sheaf $\ce_d$, for each $d$ dividing $r$, and a homomorphism 
$c_{d,d'}$ for each $d'$ dividing $d$ dividing $r$, such that 
$c_{d,d'}:\ce^{\otimes d/d'}_d \rTo \ce_{d'}$ is an isomorphism 
except where $\ce_d$ is  not locally free.  $\ce_1$ is the sheaf 
$\ck$, and the $c_{d,d'}$ should be compatible with one another 
in the obvious way (see Definition \ref{coh-quasi}). To prevent 
certain anomalies (including a potential failure to be 
separated), we also require that the length of the cokernel of 
$c_{d,d'}$ at each singularity be $d/d'-1$.  The details of these  
nets are described in Section \ref{nets}. 

In the case that the underlying curve is smooth, or in fact, 
whenever $\ce_r$ is an invertible sheaf, then the entire 
structure is canonically determined by $\ce_r$ and the 
isomorphism $c_{r,1}: \ce^{\otimes r}_r \rTo \ck$.  Obviously, 
when $r$ is prime this collection consists only of the  data of 
$\ce_r$ and $c_{r,1}: \ce^{\otimes r}_r 
\rTo \ck$. 

When the homomorphism $c_{r,1}$ fails to be an isomorphism at 
some point $\fp$, the completion of the local ring of $X$ at 
$\fp$ is of the form $\hat{\co}_{X,\fp} \cong A=k[[x,y]]/xy$, and 
$\ce_r$ corresponds to an $A$-module $E \cong <\zeta_1, \zeta_2| 
y 
\zeta_1 
=x\zeta_2 =0>$; and $c_{r,1}$ corresponds to the homomorphism
$\zeta^r_1 \mapsto x^u$, $\zeta^r_2 \mapsto y^v$, for some $u$ 
and $v$, with $\zeta^i_1 \zeta^{r-1}_2 \mapsto 0$ if $0 <i<r$.  
The restriction on the length of the cokernel amounts to saying 
that $u+v$ must be exactly $r$. Thus if $\tilde{A}=k[x] \oplus 
k[y]$ is the normalization of $A$, the data of $E$ and $c_{r,1}$ 
over $A$ correspond exactly to the choice of a \emph{free} 
$\tilde{A}$-module $\tilde{E}$ with an isomorphism 
$\tilde{E}^{\otimes r} \irightarrow (x^u,y^v)\tilde{A}$.  So 
$\ce_r$ corresponds to a sheaf $\tilde{\ce}_r$ on the 
normalization $\tilde{X} 
\rTo^{\nu} X$ of $X$ at $ \fp$ which is locally free, and $c_{r,1}$
corresponds to an isomorphism $\tilde{\ce}_r \irightarrow \nu^*
\omega_X (-up^+-vp^-)$ where $\nu^{-1} (\fp)=\{p^+,p^-\}$.  The
choices of $u$ and $v$ which may occur are, of course, controlled in
large part by degree considerations.  We will call the (unordered)
pair $\{u,v\}$ the \emph{order} of $(\ce_r,c_{r,1})$ at $\fp$. If
$\ce_r$ is locally free at $\fp$, then the order is $\{0,0\}$.

When $u$ and $v$ are relatively prime, the sheaf $\ce_r$ and the 
homomorphism $c_{r,1}$ still completely determine the remainder 
of the structure $\{\ce_d, c_{d,d'}\}$.  In fact, for any order 
$\{u,v\}$ of $(\ce_r,c_{r,1})$ the order $\{u_d, v_d\}$ of 
$(\ce_d,c_{d,1})$ is simply $\{u,v\} \pmod d$, and the pullback 
of $\ce_d$ to $\tilde{X}$ is $\tilde{\ce}^{\otimes r/d}_r 
\otimes \co(\frac{(u-u_d)}{d} p^+ + \frac{(v-v_d)}{d}
 p^-)$, with the pullback of $c_{r,d}$ being the obvious inclusion 
 homomorphism.  When $d$ does not divide $u$ (and hence $v$), the 
 order $\{u_d,v_d\}$ is not $\{0,0\}$ and 
 $\ce_d$ is not locally free.  In this case $(\ce_d, c_{d,1})$ 
 and $c_{r,d}$ are completely determined by $(\ce_r, c_{r,1})$.
 But when $\gcd (u,v)=\ell>1$  the sheaf $\ce_\ell$ is determined 
only up to its restriction to the normalization $\tilde{X}$ of 
$X$ at $\fp$.  An additional gluing datum is necessary (if $\fp$ 
is a non-separating node of $X$) to construct $\ce_\ell$ from 
$\ce_r$.  

Even when the sheaves $\ce_d$ are not locally free, to construct 
a cohomological field theory \cite{jkv,witten:r-spin-conj}, and 
in particular to have a composition axiom, as in quantum 
cohomology, we need not only roots of $\ck$, but also roots of 
$\ck \otimes \co (-\sum m_i p_i)$, for marked points $p_i$ of $X$ 
and integers $m_i$.  Letting $\bm 
= (m_1, \dots, m_n)$ be the vector of the orders of vanishing at 
the $p_i$'s, we call these {\em roots} (or {\em nets of roots}) 
{\em of $\ck$ of type $\bm$}.  We denote the stack of smooth 
roots of $\ck$ of type $\bm$ by $\stacksgrnm (\ck)$, its 
compactification using coherent nets of roots by $\stacksgrnmbar 
(\ck)$.   When $\ck = \omega$ this is the stack of $r$-spin 
curves, and is denoted $\stacksgrnmbar$.  It can be shown 
\cite{abramovich-jarvis} that the structure of such a coherent 
net of roots is 
 actually equivalent to the data of a 
 single invertible sheaf $\cl$ with an isomorphism $\cl^{\otimes r}  \rTo  
 \ck$ on a twisted (stacky) curve $\cc$ in the sense of 
 Abramovich and Vistoli 
 \cite{abramovich-vistoli:complete-moduli}, and the type $\bm$ 
 corresponds to the indices of stackiness at the marked points 
 of $\cc$.
 
 Unfortunately, the most general families of such nets of roots, without further
 restrictions, do not form a smooth stack.  The additional conditions
 we must place on them, although technical to state, amount simply to
 insisting that only the most natural and best-behaved families will
 be included in the stack.  In Section \ref{background} we recall some
 results from \cite{jarvis:spin} necessary to give these
 conditions, and in Section \ref{power} we describe these additional
 conditions, which are
 essentially equivalent to Abramovich and Vistoli's condition that
 stable maps should be \emph{balanced}.
 
 Once the definitions are in place, it is relatively straightforward
 to describe the universal deformation and thus show that this stack
 is smooth.  Moreover, there is a natural morphism from this stack
 $\stacksgrnmbar$ to the Deligne-Mumford stacks $\overline{\mbox{\sc
     Spin}}_{r,g}$ and $\mbox{\sc pure}_{r,g}$ of \cite{jarvis:spin},
 making $\stacksgrnmbar$ into the normalization of
 $\overline{\mbox{\sc Spin}}_{r,g}$ and $\mbox{\sc pure}_{r,g}$, and
 in particular, $\stacksgrnmbar$ (and$\stacksgrnmbar (\ck)$) is a smooth Deligne-Mumford stack,
 finite over $\stackmgnbar$.  This is described in Section \ref{stack}
 
 \subsubsection{Geometry}
 In Section~\ref{moduli} we treat the geometry of the stack of
 spin curves and its coarse moduli space.  In \ref{proj} we cover basic properties such as
 projectivity. In \ref{13} we discuss relations between the different
 spaces.  Finally, in \ref{14} we treat the question of
 irreducibility.  In particular, we first restrict to the complex
 numbers and use results of Mess \cite{mess:torelli} and Powell
 \cite{powell:torelli} on the Torelli group, and results of Birman
 \cite{birman:blmcg} on the kernel of the natural homomorphism from
 the mapping class group $\Gamma_{g,n}$ of genus $g$ surfaces with $n$
 punctures to the group $\Gamma_{g,n-1}$ to give a set of generators
 (finite when $g \neq 2$) for the Torelli group $\mathcal{I}_{g,n}$
 for all $g$ and $n$.  This result, and a characterization due to Sipe
 \cite{sipe:roots} of $r$-spin structures on a given curve as
 splittings of the Gysin sequence for the punctured tangent bundle,
 allows us to study the monodromy action on the set of $r$-spin
 structures of a given curve.  When $g$ is greater than one, and $r$
 (or any of the $m_i$) is odd, this action is transitive, and hence
 the moduli of $r$-spin curves of type $\bm$ is irreducible.  It is
 well-known (see \cite{mumford:theta-chars} or \cite{cornalba:theta})
 that $\overline{\scheme{S}}^{1/2}_g$ is the disjoint union of two irreducible components.
 And when $r$ and all of the $m_i$ are even, there is a natural
 morphism $[r/2]: \sgrnmbar \rTo
 \overline{\scheme{S}}^{1/2,\mathbf{0}}_{g,n}$ described in
 Section~\ref{13}, which shows that $\sgrnmbar$ must be the disjoint
 union of two at least two pieces.  But we show that when $g\geq 2$
 monodromy acts transitively on the fibres of the morphism $[r/2]$, so 
 $\sgrnmbar$ is actually the disjoint union of exactly two irreducible
 components.
 
 The case of $g=1$ has special arithmetic interest.  When
 $\bm=\mathbf{0}$ it is a classical fact that
 $\overline{\scheme{S}}^{1/r,\mathbf{0}}_{1,1}$ is the disjoint union
 of $d(r)$ irreducible components, where $d(r)$ is the number of
 divisors of $r$.  Again using our description of the Torelli group,
 we generalize this fact to show that
 $\overline{\scheme{S}}^{1/r,\bm}_{1,n}$ is the disjoint union of
 $d_{1,r}(\bm)$ components, where $d_{1,r}(\bm)$ is the number of
 divisors of $\gcd(r, m_1, \dots, m_n)$.
 
 %The irreducibility results  of Section~\ref{14} are actually
 %proved for spin curves over a general, $n$-pointed curve.  This will
 %allow us to conclude later that the boundary divisors defined in
 %Section~\ref{bdry-div} are irreducible.
 
  These results are then used to show
 that the irreducibility properties hold in any characteristic prime to $r$. 
These irreducibility properties are important in the study of the 
Picard group of $\stacksgrbar$ \cite{jarvis:picard}, and also in 
the study of the cohomological field theory related to Witten's 
conjecture \cite{jkv,witten:r-spin-conj}. 
 
%In \cite{abramovich-jarvis} and \cite{jarvis:spin} several 
%different compactifications  of the moduli of $r$-spin curves are 
%constructed, and shown to be algebraic.  The best-behaved of 
%these compactifications is the moduli of what are called {\em 
%coherent root nets on stable curves} in 
%\cite{abramovich-jarvis}, but which will be called {\em stable 
%$r$-spin curves} or just {\em 
%   $r$-spin curves} in this article.  In the case of prime $r$, 
%   these are the same as the \emph{pure spin curves} of 
%   \cite{jarvis:spin}.

 \subsection{Conventions and Notation}
 
   By a {\em curve} we mean a reduced,
 complete, connected, one-dimensional scheme over a field.  A {\em
   semi-stable curve of genus $g$} is a curve with only ordinary
 double points such that $H^1(X,\co_X)$ has dimension $g$.  And an
 {\em $n$-pointed stable curve} is a semi-stable curve $X$ together
 with an ordered $n$-tuple of non-singular points $(p_1,\dots,p_n)$,
 such that at least three marked points or double points of $X$ lie on
 every smooth irreducible component of genus $0$, and at least one
 marked point or double point of $X$ lies on every smooth component of
 genus one.  A stable (or semi-stable) curve is a flat,
 proper morphism $ X \rightarrow T $ whose geometric fibres $X_t$
 are (semi) stable curves. Except where otherwise indicated, $r$ 
 will be a fixed positive integer, and both $g$ and $n$ will be non-negative 
 integers such that $2g-2+n>0$.
 
 By {\em line bundle} we mean an invertible (locally free of rank one)
 coherent sheaf.  By \emph{canonical sheaf} we mean the relative
 dualizing sheaf of a family of curves $f:X \rightarrow T$, and this
 sheaf will be denoted $\omega_{X/T}$, or $\omega_f$, or just $\omega$.  Note that for a
 semi-stable curve, the canonical sheaf is a line bundle.  When $T$ is
 $\spec k$, for an algebraically closed field $k$, we will also write $\omega_X$ for
 $\omega_{X/T}$.

\section{Coherent Nets of Roots}\label{nets}

\subsection{Quasi-roots and Coherent Nets of 
Quasi-roots}\label{quasi}

 To begin we need the definition of torsion-free sheaves.
 \begin{defn}
  A {\em relatively torsion-free sheaf} (or just torsion-free sheaf)
  on a family of stable or semi-stable curves $f: X \rightarrow T$
  is a coherent $\co_{X}$-module $\ce$ that is flat over $T$, such
  that on each fibre $X_t =X \times_T \spec k(t)$ the induced
  $\ce_t$ has no associated primes of height one.
 \end{defn}

We will only be concerned with rank-one torsion-free sheaves.  
Such sheaves are called {\em admissible} by Alexeev 
\cite{alexeev:compact-jac} and {\em sheaves of pure dimension 
$1$} by Simpson \cite{simpson:rep-fund-group}.  Of course, on the 
open set where $f$ is smooth, a torsion-free sheaf is locally 
free. 

\begin{defn}
Given an $n$-pointed nodal curve $X/T$, a choice of $\mathbf{m} = 
(m_1,\dots,m_n)$, and a rank-one, torsion-free sheaf $\cf$ on 
$X$, an \emph{$r$th  quasi-root of $\cf$ of type $\mathbf{m}$} is 
a pair $(\ce,b)$ consisting of a rank-one, torsion-free sheaf 
$\ce$ on $X$, and an $\co_{X}$-module homomorphism $c: 
\ce^{\otimes r} \rTo \cf (-\sum m_i D_i)$ with the following properties:
\begin{enumerate}
\item $r\cdot \deg \ce= \deg \cf - \sum m_i$;
\item $c$ is an isomorphism on the open subset of $X$ where 
$\ce$ is locally free; and 
\item for each point $\fp$ of each fibre $X_t$ of $X$ where $\ce$ is not locally 
free, the length of the cokernel of $c$ at $\fp$ is $r-1$.
\end{enumerate} 
\end{defn} 

This definition of a quasi-root reduces to that of 
\cite{jarvis:spin} when the target $\cf$ is locally free.   Of 
course, when $X$ is smooth, any quasi-root is locally free and 
$c$ is simply an isomorphism.  

As explained in \cite{jarvis:spin}, the last condition is 
necessary to ensure separatedness of the stack of quasi-roots of 
a given sheaf $\cf$, and is actually a very natural condition.  
Indeed, as indicated in the introduction, on a nodal curve $X$ 
over an algebraically closed field, each singularity (point of 
$X$ where $\ce$ is not locally free) of a quasi-root $(\ce,c)$ 
uniquely determines two positive integers $u$ and $v$, summing to 
$r$:  these are the order of vanishing of $c$ on the pullback of 
$\ce$ (mod torsion) to each of the two branches of $X$ through 
the singularity.  In particular, the completion of the local ring 
of $X$ at the singularity is of the form $\hat{\co}_{X,\fp} 
\cong A=k[[x,y]]/xy$, and $\ce$ corresponds to an $A$-module $E \cong 
<\xi_1, \xi_2 |y\xi_1=x\xi_2=0>$ \cite[chap.11, 
Prop.3]{seshadri:fibres}.  If $\cf$ is locally free, it 
corresponds to the free $A$-module $A$,  and $c$ corresponds to 
the homomorphism $\xi^r_1 
\mapsto x^u$, $\xi^r_2 
\mapsto y^v$, and $\xi^i_1 \xi^{r-i}_2 \mapsto 0$ if $0<i<r$ \cite[Prop. 3.3.1]{jarvis:spin}.  If $\cf$ 
is not locally free, then it corresponds to an $A$-module $F\cong 
<\zeta_1, \zeta_2| y \zeta_1=x\zeta_2=0>$ and $c$ corresponds to 
the homomorphism $\zeta^r_1 \mapsto x^u \zeta_1,$ $ \zeta^r_2 
\mapsto y^v \zeta_2$, and $\zeta^i_1 \zeta^{r-i}_2 \mapsto 0 $ for $0 <i<r$.  
Condition (3) on the cokernel in the definition of a quasi-root 
is simply the condition that $u+v=r$ when $\cf$ is locally free, 
or $u+v=r-1$ if $\cf$ is not locally free. 

\begin{defn}\label{pair}
We will call the pair $\{u,v\}$ the \emph{order} of $(E,c)$ at 
the singularity. 
\end{defn} 

In order to produce a smooth stack, we must consider not just 
quasi-roots, but also coherent nets of quasi-roots; and we must 
place some additional conditions on the families of these nets 
which we will explain in the next section. 

\begin{defn}\label{coh-quasi}
A \emph{type-$\mathbf{m}$, coherent net of $r$-quasi-roots} of a 
rank-one torsion-free sheaf $\cf$ is a collection 
$\{\ce_d,c_{d,d'}\}$, consisting of a   rank-one, torsion-free 
sheaf $\ce_d$ for every divisor $d$ of $r$, and an $\co_X$-module 
homomorphism $c_{d,d'}:\ce^{\otimes d/d'}_{d} \rTo \ce_{d'}$ for 
each $d'$ dividing $d$ with the following properties: 
\begin{enumerate}
\item $\ce_1=\cf$ and $c_{1,1} = \mathbf{1}$
\item \label{m-vect}
For each divisor $d$ of $r$, and for each divisor $d'$ of $d$ let 
$\mathbf{m'}$ be the $n$-tuple $(m_1',\dots,m_n')$ such that 
$m'_i$ is the  unique, non-negative integer, less than $d/d'$, 
and congruent to $m_i \mod (d/d')$.  The homomorphism $c_{d,d'}$ 
is required to make $(\ce_d,c_{d,d'})$ into a $d/d'$ quasi-root 
of $\ce_{d'}$ of type $\mathbf{m'}$. 

\item The homomorphisms $\{c_{d,d'}\}$ are compatible:

For any $d''$ dividing $d'$ dividing $d$ dividing $r$, the 
homomorphism $(c_{d',d''})^{\otimes d'/d''}$ commutes with 
$c_{d,d''}$ and $c_{d',d''}$; that is the following diagram 
commutes.$$ 
\begin{diagram}
\ce^{\otimes d/d''}_d & \rTo^{c^{\otimes d'/d''}_{d,d'}} & \ce^{\otimes 
d'}_{d'} \\ &  \rdTo_{c_{d,d''}}& \dTo_{c_{d',d''}}\\ & & 
\ce_{d''}      
\\ 
\end{diagram} $$
\end{enumerate}
\end{defn}

Of course, if $\ce_r$ is locally free, then the entire net is 
completely determined by $(\ce_r,c_{r,1})$. 

\begin{nb}
%First,  we want a canonical morphism between the moduli spaces 
%$\sgrnkbar$ and $\sgsnkbar$, defined by forgetting all of the 
%terms in $\{\ce_d, c_{d,d'}\}$ except the terms involving 
%divisors of $s$.  

%The restriction we placed on $\sgsnkbar$ in 
%Definition~\ref{def-smooth-stack} required that the integers 
%$m_i$ all lie between $0$ and $s$.  As noted there, this is not a 
%real restriction, since roots of one type are uniquely determined 
%by any roots of a type congruent to them, modulo $s$.   In fact, 
%we would get an isomorphic stack by replacing the condition in 
%Definition~\ref{coh-quasi}(\ref{m-vect}) by the condition that 
%the roots all be of type $\mathbf{m}$. 

%The second, and less arbitrary, reason that we  make the 
%stipulation of Definition~\ref{coh-quasi}(\ref{m-vect})  is the 
%fact that 
  Since the type $\bm$ changes with different $d$'s dividing $r$, a
  net of roots of $\cf$ of type $\bm$ is not the same as a net of
  roots of $\cf \otimes \co(-\sum m_i p_i)$ of type $\mathbf{0}$.
  Still, for any $d$th root $(\ce_d,c_d)$ of $\cf$ of type $\bm$,
  there is a uniquely determined $d$th root $(\ce'_d,c'_d)$ of $\cf$
  of type $\bm'$ for every $\bm'$ congruent to $\bm \mod d$; namely,
  $\ce'_d:=\ce_d \otimes \co(1/d \sum (m_i-m'_i)p_i)$.  Thus the condition in
  Definition~\ref{coh-quasi}(\ref{m-vect}), requiring that
  $\mathbf{m}$ change with the different divisors is not a
  restriction.  Indeed, we would get an isomorphic stack by insisting
  that all $c_{d,d'}$ be roots of some other type congruent to $\bm
  \mod d/d'$.  One advantage to the choice we have given here is that
  it allows the stack of curves with roots to obey something like the
  composition rules of quantum cohomology.  In particular, near any
  marked point the structure should behave like a structure obtained
  by normalizing a node.  In our case, normalizing a non-locally-free
  $d$th quasi root at a singularity of the root yields a $d$th quasi
  root, whose type $m^+$ and $m^-$ at the normalized points is $u_d-1$
  and $v_d-1$ where $\{u_d,v_d\}$ is the order of the $d$th root at
  the node, and thus the type of the normalized root is bounded
  between $0$ and $d-2$, inclusive.
  When the root is locally free, composition rules like those of
  quantum cohomology still exist for these root nets, as explained in
  \cite{jkv}, but they are more subtle.

\end{nb}

Even considering coherent nets of quasi-roots is not quite 
sufficient to ensure a smooth stack.  To describe the additional 
structure necessary to make the stack smooth we first recall some  
results on the local structure of a rank-one, torsion-free sheaf 
and some results from \cite{jarvis:spin} on the local structure 
of roots. 

\subsection{Background on Local Structures} \label{background}
\subsubsection{Local Structure of Torsion-free Sheaves}

Let $R$ be the Henselization of a local ring, of finite type over 
a  field or an excellent Dedekind domain, let $\maxid$ be the 
maximal ideal of $R$.  Every nodal curve $X/\spec R$ has, at a 
node, a local ring whose Henselization is isomorphic to $A$, the 
Henselization of $R[x,y]/xy-\pi$ at $\maxid +(x,y)$ for some $\pi 
\in \maxid$.  

\begin{defn}
Over the ring $A$, for each pair $p,q \in R$, such that $pq=\pi$, 
define $E(p,q)$ to be the $A$-module generated by two elements 
$\xi_1$ and $\xi_2$, with the relations $x\xi_2=p\xi_1$, and 
$y\xi_1=q\xi_2$. 
\end{defn} 

\begin{theorem}[Faltings\cite{faltings:torfree}]
\begin{enumerate}
\item 
Any torsion-free $E$ of rank one over $A$ is isomorphic to an 
$E(p,q)$ for some $p,q \in R$ with $pq=\pi$. 
\item 
If $p,q,p',q'$ are all in $\maxid$, and $pq=p'q'=\pi$, then 
$E(p,q)$ is isomorphic to $E(p',q')$ if and only if there exists 
a unit $\alpha 
\in R^{\cross}$ such that $p'=\alpha p$, and $q'=\alpha^{-1}q$. 
\end{enumerate}  
\end{theorem} 

The condition that $E$ is an $r$th quasi-root of $A$ implies 
\cite[5.4.11]{jarvis:spin} that in the above theorem $p$ and $q$ 
can be assumed to have the property $p^u=wq^v$ in $R$, for 
$\{u,v\}$ the order of the root, as described in Definition 
\ref{pair}, and with $w 
\in R^{\times}$.  Moreover, if $w$ has 
an $r$th  root in $R^{\cross}$, then $w$ may be assumed to be 
$1$.  In particular, there is an \'etale cover of the base on 
which $w$ may be assumed to be $1$.  

\subsubsection{Local Coordinates}

In order to use Faltings' result,  we need the definition of a 
local coordinate, which is really just a way of choosing 
parameters $x$, $y$, and $\pi$. From the deformation theory of 
stable curves, we know that near a singularity $\fp$ of $X/T$ the 
complete local ring $\hat{\co}_{X,\fp}$ over $\hat{\co}_{T,t}$ is 
of the form $\hat{\co}_{X,\fp} \cong 
\hat{\co}_{T,t} [[x,y]]/(xy-\pi)$ for some $\pi \in 
\hat{\co}_{T,t}$. And over some \'etale neighborhood $T'$ of $t$, 
there is an \'etale neighborhood $U$ of $\fp$ in $X 
\cross_T T'$ with sections $x$ and $y$ in $\co_{U}$ such that
\begin{enumerate}
\item $xy=\pi \in \co_{T',t}$.
\item The ideal generated by $x$ and $y$ has the discriminant 
locus of $X/T$ as its associated closed subscheme.
\item The obvious homomorphism $\big(\co_{T',t} 
[x,y]/(xy-\pi)\big) \rightarrow
\co_{U,\fp}$ induces an isomorphism on the completions 
$\left(\hat{\co}_{T',t} [[x,y]]/(xy-\pi)\right) \irightarrow
\hat{\co}_{U,\fp}$.
\end{enumerate}

\begin{defn}\label{local}
We call such a system a {\em local coordinate} for $X/T$ near 
$\fp$. 
\end{defn}

Note that a local coordinate is not uniquely determined.  It is only
determined up to the equivalence relation generated by the following
operations:
\begin{enumerate}
\item Pullback to \'etale covers.
\item Change by units: namely $x' = \alpha x, y'=\beta y, {\pi}'
= \sigma \pi$ with $\alpha, \beta \in \co^*_{{X}'}$, and
$\alpha \beta = \sigma \in \co^*_{T'}$.
\item Switching branches; namely, interchanging $x$ and $y$.
\end{enumerate}

\subsubsection{Local Structure of $r$th Roots and the Stacks
  {\sc spin}, {\sc pure}, and {\sc root}} 

If $E=<\zeta_1, \zeta_2|p \zeta_1=x \zeta_2, y \zeta_1 =q 
\zeta_2>$ is a rank-one, torsion-free $A$-module, then any 
$A$-module homomorphism from $E^{\otimes r}$ to a torsion-free 
$A$-module must necessarily factor through the symmetric product 
$\sym^r E$ \cite[\S3.3]{jarvis:spin}.  Thus such a homomorphism 
can be defined in terms of where it takes the elements 
$$\delta_i:=\varepsilon_{r-i}:=\xi^{r-i}_1 \xi^i_2.$$

In \cite[\S5.4.1]{jarvis:spin} it is shown that for any quasi-root
$E^{\otimes r} \rTo^b A$ there exists an invertible element $a \in
A^{\cross}$ such that for any $0 \leq i \leq u$ the homomorphism $b$
takes $\delta_i$ to $a x^{u-i}p^i + \gamma_i$ and for $0 \leq j \leq
v$ the homomorphism $b$ takes $\varepsilon_i$ to $a y^{v-j} q^j +
\gamma^j$, where $\gamma_i$ and $\gamma^j$ are nilpotent elements of
$A$, and are annihilated by $\pi^r$.

``Good'' quasi-roots are those for which $\gamma_i$ and 
$\gamma^j$ are all zero, and even better are those for which 
there is some $t 
\in R$ such that $p =t^v$ and $q=t^u$; or $p=t^{v/c}$ and $q=t^{u/c}$ for $c = \gcd(u,v)$.
The former are called simply {\em roots} in \cite{jarvis:spin} 
and the latter are called {\em pure roots}; this usage differs 
from that in this paper. 
 
The stacks of $r$th roots with these conditions  form the 
compactifications $\overline{\mbox{\sc spin}}_{r,g}$ and 
$\mbox{\sc pure}_{r,g}$, respectively, of the stack of smooth 
spin curves with an $r$th root $(\ce,b)$ of $\omega$.  But 
neither of these stacks is smooth in general.  In particular, if 
$\spec\mathfrak{o}[[s_1,\dots, s_{3g-3+n}]]$ is the universal 
deformation of a curve $X$, then the universal deformation of a 
geometric point of $\overline{\mbox{\scshape Spin}}_{r,g}$ with 
underlying curve $X$ and with order $\{u_i,v_i\}$ at the node 
defined by the vanishing of $s_i$ is 
$$\spec \mathfrak{o}[[P_1,Q_1,\dots, P_k, Q_k, s_{k+1}, \dots, s_n,
s_{n+1}, \dots, s_{3g-3+n}]]/(P^{u_i}_i-Q^{v_i}_i)$$ where the 
root $\ce$ is locally of the form $E(P_i,Q_i)$, and the forgetful 
morphism $\overline{\mbox{\sc spin}}_{r,g} \rTo \mgbar$ is given 
by $s_i=P_i Q_i$.  When $\gcd (u_i,v_i)=1$ for every $i$, the 
universal deformation of a geometric point of $\mbox{\sc 
pure}_{r,g}$ is $\spec 
\mathfrak{o}[[t_1, \dots, t_k, s_{k+1}, \dots, s_{3g-3+n}]]$ where
$P_i=t^{v_i}_i$ and $Q_i=t^{u_i}_i$.  And $\mbox{\sc pure}_{r,g}$ 
is clearly the normalization of $\overline{\mbox{\sc 
spin}}_{r,g}$ in this case.  Indeed, $\mbox{\sc pure}_{r,g}$ is 
the normalization of $\overline{\mbox{\sc spin}}_{r,g}$ whenever 
$r$ is prime.  When $(u_i,v_i)>1$ this no longer holds (See Note 
\ref{gluenote})---a fact not made clear in \cite{jarvis:spin}. 

The most important properties we need to know about 
$\overline{{\mbox{\sc 
    spin}}}_{r,g}$ are that it is a Deligne-Mumford stack, such that the
forgetful morphism to $\stackmgbar$ is both proper and 
surjective, and the stack of smooth spin curves forms an open 
dense substack. Moreover, no part of the proof of these 
properties depends on special properties of $\omega$, and so for 
any line bundle $\ck$ defined on the universal curve $ \stackcgn$ 
over $\stackmgnbar$ there is a Deligne-Mumford stack 
$\mbox{\sc{Root}}^{1/r}_{g,n}(\ck)$ of curves with ``good" $r$th 
quasi-roots of $\ck$, which is proper and surjects to 
$\stackmgnbar$.  Moreover the universal deformation of a 
geometric point of ${\sc root}^{1/r}_{g,n}(\ck)$ is exactly the 
same as that for $\overline{{\mbox{\sc 
    spin}}}_{r,g}$.

\subsection{Power Maps and Roots}\label{power}

We want to describe the last condition for families of nets of 
roots necessary to make a smooth stack.  Intuitively, this  
amounts to insisting that the net be completely determined  by 
the $r$th root $(\ce_r, c_{r,d})$, whenever the intermediate 
roots are not locally free at nodes.  We will continue to use the   
notation of the previous section.
 
\subsubsection{Power Maps}
To define a general object and homomorphism of the net we will 
insist, as in the case of $\mbox{\sc pure}_{r,g}$, that there be 
a $\tau \in R$ such that $\tau^{v}=p$ and $\tau^{u}=q$ where the 
order of $\ce_r$ at the singularity defined by $\tau$ is 
$\{u\ell,v\ell\}$ for some positive integer $\ell$.  Then locally 
$\ce_r=E(\tau^{v}, 
\tau^{u})$ and $\ce_d$ should be $E(\tau^{\hat{v}}, 
\tau^{\hat{u}})$ for some $\hat{u}$ and $\hat{v}$, as explained 
below. 
 
First we give some notation.  Let $u,v$ be two non-negative integers
(not necessarily the order of $c_{r,1}$), either both zero or both
positive, and let $s=u+v$.  Given a local coordinate $(U,T',x,y,\pi)$
choose an element $\tau \in \co_{T'}$ such that $\tau^{s}= \pi$.  Now
for any non-negative integers $i$ and $j$ summing to $s$, we define
$$E_{i,j}:=E(\tau^j,\tau^i)=<\xi_1,\xi_2|\tau^j \xi_1 = x 
\xi_2,\tau^i \xi_2 =y\xi_1>.$$  Also, we let $E_{0,0}$ be the free module $$E_{0,0}:= 
<\xi_1,\xi_2|\xi_1=\xi_2>.$$  Then for any $d \geq 1$ let $u'$ 
and $v'$ be the smallest non-negative integers, congruent 
(respectively) to $d u$ and $d v \pmod {s}$, and let 
$$u''=\frac{d u-u'}{s}$$ and $$v''=\frac{d v-v'}{s}.$$ If $s$ 
divides $d u$ (and hence $d v$), then $u'=v'=0$ and $u''=du/s, 
v''=dv/s$, so $u''+v''=d$. If, however, $s$ does not divide $d 
u$, then $u'+v'=s$ and $u''+v''+1=d$.  If $u=v=0$, then 
$u'=v'=u''=v''=0$. 
 
\begin{defn}
Given a local coordinate $(U,T',x,y,\pi)$, two non-negative 
integers $u$ and $v$, either both zero or both positive, and an 
element $\tau\in 
\co_{T'}$, 
such that $\tau^{u+v} = \pi$, as above, define the \emph{$d$th  
power of $E_{u,v}$} to be the homomorphism $\phi_d:E^{\otimes 
d}_{u,v} \rTo E_{u',v'} $ as follows.

Again we let $\delta_i:=\varepsilon_{d-i}= \zeta^{d-i}_1 
\zeta^i_2$ in $\sym^d E_{u,v}$.  The $d$th  power of $E_{u,v}$ is 
the  composition of the canonical map $E^{\otimes d}_{u,v} \rTo 
\sym^d E_{u,v}$ with the map  
 $$\delta_i 
\mapsto x^{u''-i} 
\tau^{iv}\zeta_1$$ for $0 \leq i \leq u''$ and $$\varepsilon_j \mapsto 
y^{v''-j}\tau^{ju} 
\zeta_2$$ for $0\leq j\leq v''$, where $\zeta_1$ and $\zeta_2$ are 
the generators of $E_{u',v'}$, as described above. 
\end{defn}

To check that this map is well defined we must check that 
$x\phi_d(\varepsilon_j)=p\phi_d(\varepsilon_{j+1})$, $y \phi_d 
(\varepsilon_j)=q\phi_d(\varepsilon_{j-1})$, $x 
\phi_d(\delta_i)=p 
\phi_d(\delta_{i-1})$, and $ y \phi_d(\delta_{i+1})=q \phi_d(\delta_i).$  
When $u'=v'=0$, we also have potentially duplicate definitions 
since $\varepsilon_{v''}=\delta_{u''}$, so we must also check 
that $\phi_d(\varepsilon_{v''})=\phi_d(\delta_{u''})$.  These 
relations are all easy to check, except in the special cases of 
$j=v''$ or $i=u''$, where we must check, for example, 
$x\phi_d(\varepsilon_{v''}) = p \phi_d(\delta_{d-v''-1})$.  If 
$u'=v'=0$ then this amounts to showing that 
$x\tau^{v''u}\zeta_2=px 
\tau^{(u''-1)v}\zeta_1$, but $p=\tau^v$, and $\zeta_1=\zeta_2$ in 
this case, so the equality follows from the fact that 
$u''v=\frac{d u v}{s}   = v''u$.  If $u'$ and $v'$  are non-zero, 
then we must show that $x\tau^{v''u}\zeta_2 = p 
\tau^{u''v}\zeta_1$.  This follows from the definition of $u''$ 
and $v''$, as well as the relation $x\zeta_2 = 
\tau^{v'} \zeta_1$; 
in particular, $x\tau^{v''u}\zeta_2 = \tau^{v'+v''u}\zeta_1 = 
\tau^{(sv'+dvu-v'u)/s} \zeta_1 = \tau^{(vv'+dvu)/s} \zeta_1 = 
\tau^{(v(s-u') + duv)s} \zeta_1=\tau^{v+(du-u')v/s} \zeta_1 = 
\tau^{v+u''v} \zeta_1 = p\tau^{u''v} \zeta_1$.  Thus we 
have a well-defined homomorphism $\phi_d: E^{\otimes d}_{u,v} 
\rTo E_{u',v'}$.  Note that when $u=v=0$ the map $\phi_d$ is just 
the canonical isomorphism $E^{\otimes d}_{0,0} 
\cong A^{\otimes d} \irightarrow A \cong E_{0,0}$ given by 
$\zeta^d \mapsto 
\zeta$.

When $\tau$ is zero, the $d$th  power map reduces to just 
$\delta_0 
\mapsto x^{u''} \zeta_1$, and $\varepsilon_0 \mapsto y^{v''} \zeta_2$, 
and all other $\delta_i$ and $\varepsilon_i$ for $0<i,j<r$ map to 
zero.  Thus the length of the cokernel of $\phi_d$ at the 
singularity of the central fibre is just $u''+v''=d-1$ when the 
target is not free (i.e. $u'$ and $v'$ not zero and $\zeta_1 
\neq \zeta_2$), and it is $u''+v''-1=d-1$ when the target is free 
($u'=v'=0$ and $\zeta_1=\zeta_2$).  Consequently, the $d$th  
power map is always a special case of a local $d$th quasi-root of 
$E_{u',v'}$ of order $\{u'',v''\}$.  In particular, if $\ce_r$ is 
locally isomorphic to $E_{i,j}$ for $i+j=s$ and $s|r$, the $r$th 
power map takes $E_{i,j}$ to $E_{0,0}=A$ and has order $\{ri/s, 
rj/s\}$. 

\begin{nb}
Given any two non-negative integers $i$ and $j$, and any positive 
integer $d$, it is easy to see that $E_{i,j}$ has a $d$th  root  
exactly when the equivalences $d x  \equiv  i \mod {(i+j)} $ and 
$d y 
\equiv  j \mod{(i+j)} $ can be solved for positive 
integers $x$ and $y$.  And if $d$ and $i+j$ are relatively prime, 
then a unique root always exists. 
\end{nb}

\subsubsection{Roots and Nets of Roots}
\begin{defn}
{\em An $r$th  root} of a rank-one, torsion-free sheaf $\cf$ on 
$X/T$ is an $r$th  quasi-root $(\ce,b)$ of $\cf$ such that for 
each singularity $\fp$ of any fibre $X_t/t$ of $X/T$, if $\ce$ is 
not locally free at $\fp$, then there is a local coordinate 
$(U,T',x,y,\pi)$, and a choice of $\tau$ in  $\co_{T'}$, such 
that if $\{\tilde{u},\tilde{v}\}$ is the order of $b$ at $\fp$ 
and if $\ell=\gcd(\tilde{u},\tilde{v})$, then letting 
${u}=\tilde{u}/\ell$ and ${v}=\tilde{v}/\ell$, we have that $\ce$ 
restricted to $U$ is isomorphic to $E_{{u},{v}}$ and $\cf$ 
restricted to $U$ is isomorphic to $E_{u',v'}$, and such that the 
induced homomorphism $\tilde{b}:E^{\otimes r}_{u,v} 
\rTo E_{u',v'}$ is an $r$th  power map with respect to $\tau$. 
\end{defn}

Of course if $T$ is the spectrum of a field, then $\pi 
=\tau=0$ and every quasi-root is a root.  Thus the difference 
between roots and quasi-roots is in the way they vary in 
families. 

It is also worth noting that the stack of stable curves with 
$r$th roots of $\omega$ is exactly the stack $\purerg$ of 
\cite{jarvis:spin}.

\begin{defn}\label{coherent}
A {\em coherent net of roots} for a torsion-free sheaf $\cf$ on a 
nodal curve $X/T$ is a coherent net of quasi-roots 
$\{\ce_d,c_{d,d'}\}$ with the property that for each singularity 
$\fp$ of any fibre of $X/T$, if $\ce_r$ is not locally free, then 
there is a  local coordinate $(U,T',x,y,\pi)$ and a choice of 
$\tau$ in $\co_{T'}$, such that if $\{\tilde{u},\tilde{v}\}$ is 
the order of $ \ce_r$ at $\fp$ and 
$\ell=\gcd(\tilde{u},\tilde{v})$ then $\ce_r$ is isomorphic to 
$E(\tau^{v},\tau^{u})$, with $u=\tilde{u}/\ell,$ 
$v=\tilde{v}/\ell$ and every $c_{d,d'}$ is a $d/d'$-power  map 
with respect to the same parameter $\tau$. 
\end{defn}

Again, in the special case that $T$ is the spectrum of a field, 
every coherent net of quasi-roots is a coherent net of roots. If 
$r$ is prime, then a coherent $r$th-root net is simply an $r$th 
root.  Moreover, if $\ce_d$ is locally free, then $\ce_d$ 
uniquely determines all $\ce_{d'}$ and all $c_{d,d'}$ (up to 
isomorphism) such that $d'|d$.

\begin{defn}
An \emph{isomorphism of two coherent root nets} $\{\ce_d, 
c_{d,d'}\}$ and $\{\ce'_d, c'_{d,d'}\}$ is a set of isomorphisms  
$\{\alpha_d:\ce_d \rTo \ce'_d\}$ which commute with the 
homomorphisms $c_{d,d'}$ and $c'_{d,d'}$. 
\end{defn}
\begin{proposition}
The properties of being an $r$th  root or of being a coherent net 
of roots is independent of the choice of local coordinate.

\end{proposition}

\begin{proof}
It suffices to show that for any $\ce^{\otimes d} \rTo^c \cf$ for 
which there exists a local coordinate $(U,T',x,y,\pi)$ and an 
element $\tau \in \co_{T'}$ with respect to which $\ce \cong 
E(\tau^v,\tau^u)$ and $\cf \cong E(\tau^{v'},\tau^{u'})$ and $c$ 
is a $d$th  power map, and given any other local coordinate 
$(\tilde{U}, 
\tilde{T}', \tilde{x}, \tilde{y},\tilde{\pi})$, then there exists an 
\'etale neighborhood $T''$ of $t$ in $T' \cross_T  
\tilde{T}'$ and an \'etale neighborhood $U''$ of $\fp$ in $U 
\cross_X \tilde{U}$ over $T''$ and a choice of  
$\tilde{\tau}$, such that with respect to the local coordinate 
$(U'',T'',\tilde{x}, 
\tilde{y},\tilde{\pi})$, the sheaf $\ce$ is isomorphic to 
$E(\tilde{\tau}^{v}, \tilde{\tau}^{u})$, the sheaf $\cf$ is 
isomorphic to $E(\tilde{\tau}^{v'}, \tilde{\tau}^{u'})$, and $c$ 
is a $d$th  power map with respect to $\tilde{\tau}$. 

After pulling back to an appropriate \'etale neighborhood $U''$ 
over some \'etale $T''\rTo T' \cross_T \tilde{T}'$, and after 
possibly switching branches (interchanging $x$ and $y$), we may 
assume that $\tilde{x}=\alpha x$ and $\tilde{y} 
=\beta y$ for units $\alpha,\beta \in \co_{U''}$, and $\tilde{\pi} 
= \sigma 
\pi$ with $\alpha \beta = \sigma \in \co^*_{T''}$. 

If $u+v=s$, then let $a$ and $b$ in $\co^*_{U''}$ be $s$th  roots 
of $\alpha$ and $\beta$ respectively.  Then $E(\tau^v, \tau^u)$, 
with respect  to the coordinate $(x,y,\pi)$ is $<\xi_1,\xi_2| x 
\xi_2=\tau^v 
\xi_1,y \xi_2 = \tau^u \xi_2>$, which is isomorphic to 
$<\tilde{\xi}_1, \tilde{\xi}_2 | \tilde{x}\tilde{\xi}_2 = 
(ab\tau)^v \tilde{\xi}_1, \tilde{y} 
\tilde{\xi}_1 = (ab\tau)^u \tilde{\xi}_2>$, where $\tilde{\xi} = 
a^u \xi_1$ and $\tilde{\xi}_2 = b^v \xi_2$.  Thus we may take 
$\tilde{\tau}$ to be $ab\tau$, and $\ce$ is isomorphic to 
$E(\tilde{\tau}^v, \tilde{\tau}^u)$ with respect to the new 
coordinates $(\tilde{x}, \tilde{y}, \tilde{\pi})$.

Similarly, $\cf$, which is isomorphic to $E(\tau^{v'}, 
\tau^{u'})$ with respect to the old coordinates is $<\zeta_1, 
\zeta_2| x\zeta_2= \tau^{v'}\zeta_1, y \zeta_1 = \tau^{u'} 
\zeta_2>$ which is isomorphic to $<\tilde{\zeta}_1, 
\tilde{\zeta}_2 | 
\tilde{x} \tilde{\zeta}_2 = \tilde{\tau}^{v'} \tilde{\zeta}_1,
\tilde{y} \tilde{\zeta}_1 = \tilde{\tau}^{u'} 
\tilde{\zeta}_2>$, where $\tilde{\zeta}_1 = a^{u'} \zeta_1$, and 
$\tilde{\zeta}_2  =b^{v'} \zeta_2$.  Thus, $\cf$ is isomorphic to 
$E(\tilde{\tau}^{v'}, \tilde{\tau}^{u'})$ with respect to the new 
coordinates $(\tilde{x}, \tilde{y}, \tilde{\pi})$.  Moreover, the 
$d$th  power homomorphism $\phi_d$ with respect to the old 
coordinates gives $\delta_i \mapsto x^{u''-1} \tau^{iv} \zeta_1  
= \tilde{x}^{u''-i} \tilde{\tau}^{iv}\tilde{\zeta}_1 
\alpha^{i-u''} (a b)^{-iv} a^{-u'}$; whereas the $d$th  power 
homomorphism $\psi_d$ with respect to the new coordinate gives 
$\delta_i = b^{-iv} a^{-(d-i)u} \tilde{\delta}_i \mapsto b^{-iv} 
a^{(i-d)n} \tilde{x}^{u''-i} \tilde{\tau}^{iv} \tilde{\zeta}_1$.   
These agree because $u+v=s$ and $su''+u'=du$.  A similar 
computation shows that the $d$th  power maps agree on 
$\varepsilon_i$, and thus the diagram 
\begin{diagram}
(E(\tau^v,\tau^u)_{old})^{\otimes d} & \rTo^{\phi_d} & 
E(\tau^{v'}, \tau^{u'})_{old} \\
\dTo & & \dTo\\
(E(\tilde{\tau}^v, \tilde{\tau}^u)_{new})^{\otimes d} & 
\rTo^{\psi_d} & E(\tilde{\tau}^{v'}, \tilde{\tau}^{u'})_{new}\\ 
\end{diagram} 
commutes.
\end{proof}

%Although the conditions on roots being defined by power maps with  
%respect to some parameter seem quite technical, they are, in 
%fact, simply requiring that the coherent root net behave in 
%the most natural way possible.  The main result of this section 
%is that by requiring   this ``most natural" behavior, we actually 
%get a smooth, Deligne-Mumford stack, proper and flat over 
%$\stackmgnbar$. 

\begin{example}
When the target sheaf is $\omega$ and $r$ is two, then coherent 
nets on smooth curves correspond to classical spin curves (a 
curve and a theta-characteristic) \emph{together with an explicit 
isomorphism $\ce^{\tensor 2}_{2}  \irightarrow \omega$.} 
\end{example}

\begin{example}\label{target}
If the target sheaf is $ \co_X$, then a coherent net of $r$th  
roots of $\co$ on a smooth curve $X$ corresponds to an 
$r$-torsion point of the Jacobian of $X$, again with an explicit 
isomorphism of $\ce^{\tensor r}_{r}  
\irightarrow  
\co_X$.  In particular, the stack of 
coherent nets of $r$th  roots of $\co$ on smooth curves of genus 
one (with one marked point) forms a gerbe over the disjoint union 
$\coprod_{d|r} Y_1(d)$ of the modular curves $Y_1(d)$ of points 
of exact order $d$, where $d$ divides $r$.  Since the stack of 
nets of roots of $\co$ on stable curves is the normalization with 
respect to $\stackmgnbar$ of the stack of nets on smooth curves 
(see Theorem \ref{line}), it forms a gerbe over the union of 
modular curves $\coprod_{d|r} X_1(d)$. 

\end{example}

It is easy to see that if $\ce_d$ is locally free, then for each 
$d'$ dividing $d$, $\ce_{d'}, c_{d',1}$, and $c_{d,d'}$ are all 
canonically determined (up to isomorphism) as $\ce_{d'} = 
\ce^{\otimes d/d'}_d$, $c_{d',1} = c_{d,1}$, and $c_{d,d'}= \mathbf{1}$. 
Similarly if $\ce_{d'}$ is not locally free, then $\ce_{d'}$ and 
$c_{d,d'}$ for each $d$ such that $d'|d|r$ are canonically 
determined by $\ce_r$ and $c_{r,1}$.  In particular, off the 
singularities we have $\ce_{d'}= \ce^{\otimes 
d/d'}_d=\ce^{\otimes r/d'}_r$, and near the singularities 
$\ce_{d'}$ is determined as a $d/d'$-power of $\ce_d$.  Thus 
\'etale descent shows that $\ce_{d'}$ is globally determined (up 
to isomorphism) by $\ce_r$.  And this occurs precisely when $d'$ 
does not divide $\tilde{u}$ (or $\tilde{v}$), where 
$\{\tilde{u},\tilde{v}\}$ is the order of $c_{r,1}$. 
 
In fact, the only additional information we gain from the 
coherent net of roots that is not inherent in the root 
$\ce^{\otimes r}_r 
\rTo \cf$ is in the case that $\cf$ is locally  free 
and$(u,r)=\ell>1$.  In this case, the coherent root net amounts 
to the choice of a locally free $\ell$th  root $\ce_\ell$ of 
$\cf$ and a non-locally-free $(r/\ell)$th  root of $\ce_\ell$. 
 
The actual objects we wish to study are not just coherent root 
nets on a fixed curve, but rather, fixing an invertible sheaf 
$\ck$ on the universal curve $\stackzgn$ over $\stackmgnbar$, we 
wish to study pairs $(X/T,\{\ce_d,c_{d,d'}\})$ of a curve $X/T$ 
and a coherent net of $r$th roots $\{\ce_d, c_{c,c'}\}$ of the 
line bundle $\ck|_{X}$.  Of special interest is the case of $\ck 
= \omega_{X/T}$ and its twists by canonical sections. 

\begin{defn}
A genus-$g$, $n$-pointed $r$-spin curve of type $\bm$ is a stable 
$n$-pointed curve $X/T$ of genus $g$  and a coherent net of $r$th 
roots of $\omega(\bm):=\omega_{X/T}(-\sum m_i D_i)$, where $D_i$ 
is the divisor of $X$ corresponding to the $i$th section 
$\sigma_i:T 
\rTo X$. 
\end{defn}

\begin{defn}\label{iso}
  An isomorphism of $r$-spin curves $(X/T, \{\ce_d, c_{d,d'}\})$ and
  $(X'/T', \{\ce'_d, c_{d,d'}\})$ is a pair of isomorphisms
  $(i,\alpha)$ where $i:X/T \rTo X'/T'$ is an isomorphism of stable
  curves, and $\alpha:\{\ce_d, c_{d,d'}\} \rTo
  \{i^*\ce_d,i^*c_{d,d'}\}$ is an isomorphism of coherent nets of
  roots, compatible with the canonical isomorphism $\omega_{X/T} (\bm)
  \irightarrow i^* \omega_{X'/T'}(\bm)$.
\end{defn} 

\begin{example}
As explained in example \ref{target}, if $n \ge 1$ and $\mathbf m 
= 
\mathbf 0$ then $\omega(\bm)=\co$, and a smooth $n$-pointed 
$r$-spin curve of genus $1$ is determined  by a $n$-pointed curve 
of genus  $1$ with a point of order $r$ on the curve.  However, 
the automorphisms of the underlying curve identify some of these 
$r$-spin structures.  In particular, when $n=1$ and $r$ is odd, 
the elliptic involution acts freely on all the non-trivial 
$r$-spin structures, thus there are only $1+(r^2-1)/2$ 
isomorphism classes of $r$-spin structures on the generic 
$1$-pointed curve of genus $1$. 
\end{example}

\subsection{The Stack of Coherent Nets of Roots on 
Curves}\label{stack}
\begin{defn}
For any $g$ and $n$ such that $2g-2+n>0$, and for any $n$-tuple 
of integers $\bm=(m_1, \dots, m_n)$, if $\ck$ is a line bundle on 
the universal curve $\cgnvbt$ over $\mgnbar$, then we define 
$\sgrnmkbar$ to be the stack of genus $g$, $n$-pointed, stable 
curves, together with the data of a coherent net of $r$th-roots 
of $\ck$ of type $\bm$ on the curve.  Isomorphisms of these 
objects are analogous to those described in Definition \ref{iso} 
(where $\omega(\bm)$ is replaced by a general line bundle $\ck$). 
\end{defn}

The main result of this section is that $\sgrnmkbar$ is a smooth 
Deligne-Mumford stack over $\mathbb{Z}[1/r]$, which is proper 
(finite) over $\mgnbar$. 

We begin with some deformation theory.  In \cite{jarvis:spin} it 
is claimed that the universal deformation of a pure spin curve 
over the spectrum of a field is the $r$-fold cover of the 
universal deformation of the underlying curve 
$$\spec \mathfrak{O}[[\tau_1,\dots, 
\tau_m,t_{m+1}, \dots,  t_{3g-3}]] \rTo \spec \mathfrak{O} [[t_1, 
\dots, t_{3g-3}]]$$ given by $t_i=\tau^r_i$ when $i\leq m$ (and 
$\{\tau_i=0\}$ are the loci of singularities of the underlying 
curve where the spin structure is not locally free). 

This claim is correct for $r$ prime, and even for $r$ composite, 
provided that no singularity has order $\{u,v\}$ with $\gcd 
(u,v)>1$. However, when $\gcd (u,v)>1$, this is no longer 
correct, and the additional datum of the coherent net of roots is 
precisely what is needed to remedy the defect. 

\begin{theorem}\label{univ}
Let $\overline{X}$ be a genus-$g$, stable, $n$-pointed curve, and 
let $\ck$  be a line bundle defined on the universal universal 
deformation of the curve $\overline{X}$.  Let 
$\{\overline{\ce}_d,\overline{c}_{d,d'}\}$ be a coherent net of 
$r$th roots of the  restriction $\overline{\ck}$ to 
$\overline{X}$. The universal deformation space of the curve 
$\overline{X}$ with the net $\{\overline{\ce}_d, 
\overline{c}_{d,d'}\}$ is the cover 
 
$$\spec \mathfrak{o} [[\tau_1, \dots, \tau_m, t_{m+1} \dots, 
t_{3g-3+n}]] \rTo \spec \mathfrak{o}[[t_1, \dots, t_m, t_{m+1}, 
\dots, t_{3g-3+n}]]$$ where $t_i=\tau^{r_i}_i$, and $r_i = r/
\gcd(u_i,v_i)$, and $\{u_i,v_i\}$ is the order of the $i$th  
singularity of $(\overline{\ce}_r, \overline{c}_{r,2})$. In 
particular, for any $\ck$ defined on the universal curve over 
$\mgnbar$, the stack of curves with a coherent net of roots of 
$\ck$ is smooth. 
\end{theorem} 

\begin{proof}
  Let $R= \mathfrak{o}[[\tau_1, \dots, \tau_m, t_{m+1}, \dots,
  t_{3g-3+n}]]$.  There exists a unique curve $X$ over $\spec R$
  induced from the universal deformation of $\overline{X}$ over
  $\mathfrak{o}[[t_1, \dots t_{3g-3+n}]]$.  Moreover, on the smooth
  locus $V \subseteq X$ of $f: X \rTo \spec R$, there is a unique
  extension of the net of roots $\{\overline{\ce}_d, \overline{c}_{d,d'}\}$ on
  the special fibre to a net of roots $\{\ce_{d,V}, c_{d,d',V}\}$ on
  $V$.  Similarly, about each singularity of the special fibre, there
  is a local coordinate with $\pi= t_i$.  
  $(\overline{\ce}_r,\overline{c}_{r,1})$ uniquely
  determine the order $\{\tilde{u},\tilde{v}\}$, so one
  may define the obvious coherent root net on the neighborhood
  $U_i$ of the local coordinate; namely, if
  $\ell=\gcd(\tilde{u},\tilde{v})$ and $u = \tilde{u}/\ell$,
  $v=\tilde{v}/\ell$, and if $u'$, and $v'$ are defined as the
  smallest non-negative integers congruent $\mod (u+v)$ to $du$ and
  $dv$, respectively, then let $\ce_{d,U_i}:= E_{u',v'}$.  The
  homomorphism $c_{d,d'}$ is defined to be the obvious power map.  To glue these different nets on $V$ and $U_i$
  together requires descent data, in the form of $r$th roots of the
  transition functions $\sigma_i$ of $\ck$ on the intersections $U_i
  \cap V$.  On the special fibre these data exist; namely for each
  $d$, local isomorphisms $\overline{\ce}_d \irightarrow E(0,0)$
  uniquely determine descent data on the special fibre as
  $\overline{\gamma}_{d,i}$ which are $d$th roots of
  $\overline{\sigma}_i$.

Now $\overline{\gamma}_{d,i}$ has a unique extension to some 
$\gamma_{d,i}$ defined on an \'etale cover of $U_i \cup V$.  This  
is true since it holds on the strict Henselization of the local 
ring of any point, and since $U_i \cup   V$ is quasi-compact, the 
small 
\'etale site $(U_i \cup V)_{\acute{e}t}$ is Noetherian, and thus 
there exists an \'etale cover of $U_i \cup V$ on which 
$\gamma_{d,i}$ can be defined.  The fact that $r$ is invertible 
implies that the choice of $\gamma_{d,i}$ is unique.  Thus we 
have uniquely determined descent data for $\{\ce_{d,V}, 
c_{d,d',V}\}$ and $\{\ce_{d,U_i}, c_{d,d'U_i}\}$ and thus a 
canonically defined, coherent net $\{\ce_d, c_{d,d'}\}$ for 
$\ck$. 

For any deformation $(X,\{\ce_d,c_{d,d'}\})$ of 
$(\overline{X},\{\overline{\ce}_d,\overline{c}_{d,d'}\})$ over an 
Artin local ring $S$, there is a canonical morphism of $\spec S$ 
to $\spec 
\mathfrak{o}[[t_1, 
\dots, t_{3g-3+n}]]$ induced by the underlying curve.  And for 
each singularity, (say defined by $t_i=0)$ there is a choice of 
an element $\alpha_i 
\in S$ such that $\alpha_i^{r_i} = t_i$, such that the 
deformation is a net of power maps with respect to $\alpha_i$.  
Thus we have a homomorphism $R \rTo S$ (defined by $\tau_i 
\mapsto \alpha_i$) lifting the map $\mathfrak{o}[[t_i, \dots, 
t_{3g-3+n}]]\rTo S$. 

The induced net of roots over $X/S$ is locally isomorphic to the 
given deformation $\{\ce_d, c_{d,d'}\}$.  But the gluing data 
induced by this choice of $\alpha_i$ may not be the same as those 
induced by $\{\ce_d\}$.  Nevertheless, these gluing data differ 
at worst by an $s$th  root of unity; and, as explained below, one 
can replace $\alpha_i$ with $\rho \alpha_i$, for an appropriate 
$s$th  root  of unity $\rho$, to get a different morphism $R \rTo 
S$ which will induce not only local isomorphisms of the induced 
net with $\{\ce_d, c_{d,d'}\}$, but which will also induce global 
isomorphisms. 

This shows that the deformation over $\spec R$ is versal.  Since 
the isomorphism functor for coherent nets is unramified (Lemma
\ref{unram}), the deformation is actually universal. 
\end{proof}

\begin{nb}\label{gluenote}
  The choice of gluing data is uniquely determined once an isomorphism
  $\overline{\ce}_r \irightarrow E(0,0)$ is chosen at each singularity
  of the central fibre; however, a different isomorphism can yield
  potentially different gluing data.  The difference between two sets
  of gluing data induced in this way is root of unity; indeed, any
  automorphism of $E(0,0)$ compatible with a $d$th power map to a
  rank-one free module is easily seen to be of the form $\left(
    \begin{array}{cc} \alpha & 0 \\ 0 & \beta
\end{array} \right) $ with $\alpha^d = \beta^d =1$, (that is $E(0,0) = <\zeta_1, 
\zeta_2|x\zeta_1 = y\zeta_1=0>$ is mapped by $\zeta_1 \mapsto \alpha
\zeta_1, \zeta_2 \mapsto \beta \zeta_2)$, and such an automorphism
induces a change in $\overline{\gamma}$ on the $x$ branch by 
$\alpha$ and a change on the $y$ branch by $\beta$.  In 
particular, if the order of the singularity of $\overline{\ce}_r$ 
is $\{\tilde{u},\tilde{v}\}$ and if 
$\gcd(\tilde{u},\tilde{v})=\ell$ and $s=r/\ell$, then since 
$\overline{\ce}_r$ has an $s$th power map to 
$\overline{\ce}_\ell$, which is locally free near the 
singularity, all automorphisms must induce a change in gluing 
data which is an $s$th root of unity. 

However, if $v=\tilde{v}/\ell$ and $u=\tilde{u}/\ell$, the 
automorphism $\left( 
\begin{array}{cc} 
\alpha & 0 \\ 0 & \beta 
\end{array} \right) $ extends to an isomorphism $E(t^{v}, t^{u}) 
\irightarrow E( t^{v}\alpha/\beta, t^{u}\beta/\alpha )$.  And, 
since $\gcd(u,v)=1$, and $u+v=s$, there is a unique $s$th  root 
of unity $\rho$ such that $\rho^{v}=\alpha/\beta$ and $\rho^{u} 
= \beta/\alpha$.  Thus  $\left( \begin{array}{cc} 
\alpha & 0 \\ 0 & \beta 
\end{array} \right) $ induces an isomorphism $E(t^{v}, t^{u}) 
\irightarrow E(\tilde{t}^{v}, \tilde{t}^{u})$ where $\tilde{t} = \rho t$.  
Consequently, the different set of gluing data induced by the 
automorphism   $\left( \begin{array}{cc} 
\alpha & 0 \\ 0 & \beta 
\end{array} \right) $ could also be induced by making a different 
choice of the parameter $t$, but leaving the isomorphism 
$\overline{\ce}_r 
\irightarrow E(0,0)$ unchanged. 

This also illustrates the flaw in the argument of 
\cite{jarvis:spin} that $\purerg$ has a smooth universal 
deformation space; in that case there is no $s$th power map, so 
an automorphism may induce an $r$th root of unity which is not an 
$s$th root of unity, and thus cannot be induced by a different 
choice of $t_i$.  In particular, the deformation given there is 
not versal. 
\end{nb}

%Recall first the definition of a \emph{Deligne-Mumford morphism} 
%from \cite[3.1.2]{jarvis:pdg}. 

%\begin{defn}
%A morphism of stacks $f: F\rTo G$ is called 
%\emph{Deligne-Mumford} if for every morphism $X \rTo G$ from a 
%representable stack $X$ to $G$, the fibre product $F \cross_G X$ 
%is a Deligne-Mumford  stack.  Alternately we may say that $F$ is 
%\emph{relatively Deligne-Mumford} over $G$. 
%\end{defn}

%The basic property of Deligne-Mumford morphisms proved in 
%\cite{jarvis:pdg} is \emph{If $G$ is Deligne-Mumford and $F$ is 
%relatively Deligne-Mumford over $G$, then $F$ is 
%Deligne-Mumford.}

We can now prove the main theorem of this section.

\begin{theorem}\label{line}
  For any line bundle $\ck$ on the universal curve $\stack{C}_{g,n}$
  over $\mgnbar$, the stack $\stacksgrnmkbar$ is the normalization of 
  the stack of ``good'' quasi-roots 
  $\mbox{\sc Root}^{1/r}_{g,n} (\ck \otimes \co(-\sum p_i m_i))$ (see Section \ref{local}), and
  in particular, it is a smooth, proper Deligne-Mumford stack over
  $\mathbb{Z}[1/r]$, and the natural forgetful morphism $\sgrnkbar
  \rTo \mgnbar$ is finite and surjective.
\end{theorem} 

\begin{proof}
  The proof follows directly from the description of the universal
  deformation and the fact that the diagonal is representable, proper,
  and unramified (Lemma \ref{unram}).  In fact, it is easy to produce
  a scheme $T''$ which is \'etale over $\stacksgrnmkbar$, since there
  is a scheme $T$ which is an \'etale cover of $\mbox{\sc
    Root}^{1/r}_{g,n}(\ck \otimes \co(-\sum m_ip_i))$, and for which
  there is a local coordinate such that ``good" quasi-roots $(\ce,b)$
  of $\ck \otimes \co(-\sum m_i p_i)$ have the form $E(p,q)$ with
  $p^u=q^v$, and $b$ is defined as $b(\delta_i)=x^{u-i}p^i$ for $0
  \leq i \leq u$ and $b(\varepsilon_j)=y^{v-j}q^j$.

If $T''$ is the scheme locally defined as $\spec 
\prod_{\mu_{\ell}} 
\co_T [t]/(t^{v/\ell}-p, t^{u/\ell}-q)$ where $\{u,v\}$ is the 
order of $(\ce,b)$ at the singularity, and $\ell=\gcd(u,v)$, then 
letting the   universal net of $r$th roots defined on $T''$ be 
the obvious one, with the different choices of gluing for 
$\ce_{\ell}$ indexed by the elements of $\mathbf{\mu}_{\ell}$, it 
is clear that $T''$ is \'etale over $\stacksgrnmkbar$, and that 
$T''$ is the normalization of $T$. 
\end{proof}

\begin{lemma}\label{unram}
  
  The relative (over $\stackmgnbar$) diagonal $\Delta: \sgrnkbar
  \cross \sgrnkbar \cross_{\mgnbar} T\rTo \sgrnkbar \cross_{\mgnbar}
  T$ is representable, proper, and unramified.  That is to say, given
  a stable curve $X/T$, the functor of isomorphisms of coherent nets
  of roots of $\ck$ on $X/T$ is representable, proper, and unramified.
\end{lemma}

Of course, since the diagonal $\stackmgnbar \cross \mgnbar \rTo
\mgnbar$ is representable, proper, and unramified
\cite{deligne-mumford}, this lemma shows that the (absolute) diagonal
$\sgrnkbar \cross \sgrnkbar \rTo \sgrnkbar $ also has those
properties.

\begin{proof}
  
  The functor of isomorphisms of curves and quasi-roots of a given
  line bundle $\ck$ has all of the stated properties
  \cite[Propositions 4.1.14, 4.1.15, and
  4.1.16]{jarvis:spin}.\footnote{Again the proofs there are given only
    for quasi-roots of the line bundle $\ck = \omega$, but they do not
    depend on any property of $\omega$ except the fact that it is a
    line bundle on the universal curve, thus the results are true for
    a general $\ck$.}  But an isomorphism of a coherent net of roots
  is just a coherent system of isomorphisms of the underlying
  quasi-roots.  Thus the functor of isomorphisms of nets of roots is 
representable as the locus where the individual isomorphisms of 
the underlying quasi-roots all agree.  

To check properness we use the valuative criterion and check that 
any isomorphism on the generic fibre extends to the whole curve.  
Since this holds for the individual isomorphisms of the 
underlying quasi-roots, the only additional thing to check is 
that the individual isomorphisms of the terms in the net are all 
compatible. But it is clear in the case that the base is a 
discrete valuation ring that any isomorphisms that are compatible 
on the generic fibre must also be compatible on the entire curve. 

To check that the functor of isomorphisms is unramified, it suffices
to check that any automorphism of a coherent root net over a ring $R$
with square-zero ideal $I$, such that the automorphism is the identity
over $\overline{R} = R/I$ is the identity over $R$.  But again, the
automorphisms of the underlying quasi-roots must be trivial on $R$
when they are trivial on $\overline{R}$, and thus
the automorphism of the net is also trivial on $R$.
\end{proof}

\section{Geometry of Spin Curve Stacks and Their Moduli}\label{moduli}

Our chief interest is in $r$-spin curves, that is when $\ck$ is the
canonical (relative dualizing) bundle of the universal curve.  We
denote the stack of $n$-pointed, genus-$g$, stable $r$-spin curves of
type $\bm$ by $\stacksgrnmbar$, smooth spin curves by $\stacksgrnm$
and their coarse moduli spaces by $\sgrnmbar$ and $\sgrnm$,
respectively.  In this section we study some basic geometric
properties of these spaces and the relations between them, and we
describe the configuration of their irreducible (and connected)
components.

\subsection{Basic Properties}\label{proj} 

 Since the stacks are all Deligne-Mumford, they have 
coarse moduli spaces which are {\em a priori} only algebraic 
spaces \cite[Corollary 1.3.1]{keel-mori}. 
%fix-double-check theorem reference number
But it is straightforward to see that in our case, the moduli 
spaces are actually projective schemes. 

\begin{proposition}
The moduli spaces $\sgrnmbar$ and $\sgrnm$ are normal and 
projective (respectively, quasi-projective). 
\end{proposition}

\begin{proof}
  The coarse moduli space of any smooth algebraic stack is normal.
  (This can be seen from the proof of \cite[Proposition
  2.8]{vistoli}.)
  
  The natural forgetful map $\sgrnmbar \rightarrow \schememgnbar$ to
  the moduli $\schememgnbar$ of stable $n$-pointed curves is
  surjective and finite.  Thus any ample bundle on $\schememgnbar$
  pulls back to an ample bundle on $\sgrnmbar$ (see \cite[proof of
  3.11]{kollar:projectivity}).  But it is well known that
  $\schememgnbar$ is projective (e.g., \cite[Theorem
  6.40]{harris:moduli-book}), and thus $\sgrnmbar$ is also.

\end{proof}

 \subsection{Relations Between the Different Spaces}\label{13}

There are several natural morphisms between the various stacks 
(and moduli spaces). 
\begin{enumerate}
\item There is a canonical isomorphism from
  $\overline{\stack{S}}^{1/r,\bm}_{g,n}$ to
  $\overline{\stack{S}}^{1/r,\bm'}_{g,n}$ where $\bm'$ is $n$-tuple
  whose entries are all congruent to $\bm \mod r$; namely for any net
  $\{\ce_d, c_{d,d'}\}$ of type $\bm$, let $\{\ce'_d, c'_{d,d'}\}$ be
  the net given by $\ce'_d =\ce_d \otimes \co(1/d \sum (m_i-m'_i)p_i)$
  where $p_i$ is the $i$th marked point, and $c'_{d,d'}$ is the
  obvious homomorphism.  Because of this canonical isomorphism, we
  will often assume that all the $m_i$ lie between $0$ and $r-1$
  (inclusive).
\item The universal curve over $\stacksgrnmbar$ is easily seen to be
  the stack $\stack{S}_{g,n+1}^{1/r,\bm'} \rTo^{\pi} \stacksgrnmbar$,
  where $\bm'$ is the $(n+1)$-tuple $(m_1, \dots, m_n, 0)$; and $\pi$
  is the morphism which simply forgets the $(n+1)$st marked point.  If
  $m_{n+1}$ is not congruent to zero $\mod r$, there is no such
  morphism, since the degree of $\omega(\bm')$ is
  $2g-2-\sum^{n+1}_{i=1}m_i$ and the degree of $\omega(\bm)$ is
  $2g-2-\sum^n_1 m_i$, but both cannot be simultaneously divisible by
  $r$, and thus at least one stack is empty.
\item If $s$ divides $r$ and $\bm'$ is the $n$-tuple $(m'_1, \dots,
  m'_n)$ such that $m'_i$ is the least non-negative integer congruent
  $m_i \mod s$, then there is a natural map $\stacksgrnmbar
  \xrightarrow{[r/s]} \overline{\stack{S}}^{1/s,\bm'}_{g,n}$ defined
  by forgetting all of the terms in the net except those indexed by
  divisors of $s$.

Over a smooth curve, since the nets are determined by 
$(\ce_r,c_{r,1})$ (or $(\ce_s,c_{s,1})$ for 
$\overline{\stack{S}}_{g,n}^{1/s,\bm}$), this morphism $[r/s]$ is 
equivalent to replacing $\ce_r$ with $\ce^{\otimes r/s}_r 
\otimes \co (1/s\sum((m_i-m'_i)p_i))=\ce_s$, and $c_{r,1}$ with $c_{r,1} \otimes i$, where $i: \co \rTo \co
(\sum (m_i-m'_i)p_i)$ is the canonical inclusion.  In the case 
that $r/s$ and $s$ are relatively prime, $\stacksgrnm$ is 
actually a product 
$$\stacksgrnm \cong \stack{S}^{1/s,\bm'}_{g,n} \cross_{\mgn}
\stack{S}^{s/r,\bm''}_{g,n}.$$
This follows from the fact that on one
hand, the obvious maps $\stacksgrnm \rTo \stack{S}^{1/s,\bm'}_{g,n}$
and $\stacksgrnm \rTo \stack{S}^{s/r,\bm''}_{g,n}$ induce a map from
$\stacksgrnmbar$ to the product.  And on the other hand, the inverse
of this map can be constructed as $((X,\cl_s,c_s),(X,\cm_d,c'_d))
\mapsto (X,\cl^{\tensor a}\tensor\cm^{\tensor b},c^{ad}_s \otimes
c^{'bs}_d)$, where $ad+bs = 1$.  It is not difficult to check that
this map does not depend on the specific choice of $a$ and $b$, and so
is well-defined.
\end{enumerate}

 \subsection{Irreducibility}\label{14}
\subsubsection{Irreducibility over $\mathbb{C}$}

In the special case when $r$ is $2$ and $g \geq 1$, it is known 
(see \cite{cornalba:theta} or \cite{mumford:theta-chars}) that 
$\overline{\scheme{S}}^{1/2}_g $ is the disjoint union of two 
irreducible components $\overline{\scheme{S}}^{1/2 
\text{ even}}_g$ and $\overline{\scheme{S}}^{1/2 \text{ odd}}_g$ corresponding to the even 
and odd theta characteristics, respectively.  We now extend these 
results for general $g$, $r$, $n$, and $\mathbf{m}$, assuming the 
base field is $\mathbb{C}$.  

If $g=0$, let $\ell_{0,r}(\mathbf{m})=1$, if $g=1$, let
$\ell_{1,r}(\mathbf{m})$ be $\gcd(r,m_1, \dots, m_n)$, and if $g \neq
2$, let $\ell_{g,r} (\mathbf{m})$ be $\gcd (2,r, m_1, \dots, m_n)$;
thus $ \ell_{g,r}(\mathbf{m})$ is either $1$ or $2$ for $g \geq 2$.
For any genus $g$, let $d_{g,r} (\mathbf{m})$ be the number of
positive divisors of $\ell_{g,r} (\mathbf{m})$, including
$\ell_{g,r}(\mathbf{m})$ and $1$.

\begin{theorem}\label{disjoint}
 The moduli space $\sgrnmbar$ is the disjoint union of 
$d_{g,r}(\mathbf{m})$ irreducible components. 
\end{theorem}

In the case of $g=0$, $\omega(\bm)$ has a unique $r$th root, if 
any, and so the coarse moduli 
$\overline{\scheme{S}}^{1/r,\bm}_{g,n}$ is either empty or is 
isomorphic to the coarse moduli space 
$\overline{\scheme{M}}_{0,n}$ for all $r$, $n$, and $\mathbf{m}$.  
In either case it is irreducible.  Of course the additional 
structure of the homomorphisms in the net prevents the  stack 
$\stack{S}^{1/r, 
\bm}_{0,n}$ from being isomorphic to the stack $\overline{\stack{M}}_{0,n}$. 

Consider now the case that $g \geq 1$.  First, to see that the 
number of irreducible (and connected) components of $\sgrnmbar$ 
is at least $d_{g,r}(\mathbf{m})$, consider the morphism 
$$[r/\ell_{g,r}(\mathbf{m})]:\sgrnmbar \rTo
\overline{\scheme{S}}^{1/{\ell_{g,r}},\bm'}_{g,n}$$
Since $\ell_{g,r}(\mathbf{m})$ divides all of the $m_i$, it will 
be enough to show that $\sgrnmbar$ has $d_{g,r}(\mathbf{m})$ 
disjoint irreducible components when $r$ divides all of the 
$m_i$.  But in this case $\sgrnmbar$ is canonically isomorphic to 
$\overline{\scheme{S}}^{1/r,\mathbf{0}}_{g,n}$, so it suffices to 
assume $\mathbf{m}=\mathbf{0}$. 

In the case of $g=1$, $\omega$ is trivial, and $r$th roots of 
$\omega$ are simply $r$-torsion points of the Jacobian.  This 
case follows from a classical fact. 
\begin{lemma} \label{elliptic} The 
moduli space of elliptic curves with an $r$-torsion point 
consists of one irreducible connected component for each positive  
divisor of $r$.
\end{lemma} 

Indeed, if $\jac X$ is represented as the quotient of
$\mathbb{C}/\Lambda$ for some lattice $\Lambda=<1,\tau>$, and
$\jac_rX$ is represented as the quotient $(1/r)\Lambda/\Lambda$,
giving an isomorphism $\jac_r X \irightarrow \mathbb{Z}/r\mathbb{Z}
\cross \mathbb{Z}/r\mathbb{Z}$, then a point corresponding to $(a,b)
\in \mathbb{Z}/r\mathbb{Z} \cross \mathbb{Z}/r\mathbb{Z}$ is in the
irreducible component indexed by the divisor $\gcd(r,a,b)$
\cite[V.3.3]{farkas-kra}.

In the case that $g \geq 2$ then $\ell_{g,r}(\mathbf{m})$ is 
either $1$ or $2$.  In the first case there is nothing to prove, 
and in the second case the space 
$\overline{\scheme{S}}^{1/\ell_{g,r},\mathbf{0}}_{g,n}$ is 
$\overline{\scheme{S}}^{1/2,\mathbf{0}}_{g,n}$.  In this  case we 
have the following lemma, which is a generalization of its 
classical counterpart, and which completes the proof that for all 
$g$, $r$, $n$,  and $\bm$ the number of irreducible and connected 
components is at least $d_{g,r}(\mathbf{m})$. 

\begin{lemma} \label{def-invar}
If $m_i$ is even for every $i$, then the function 
$e:\overline{\scheme{S}}^{1/2,\bm}_{g,n} 
\rightarrow \mathbb Z/2\mathbb Z$, which takes $ (X,(\ce_2,b))$ to $ \dim
H^0(X,\ce_2) \pmod 2$, is locally constant.
\end{lemma}

The proof of Lemma~\ref{def-invar} is an easy generalization of 
Mumford's proof \cite[Section 1]{mumford:theta-chars} of the 
corresponding result when $n=0$.  Mumford's idea is to take a 
divisor $D$ of high degree and make a quadratic form $q$ on 
$H^0(X,\ce_2(D)/\ce_2(-D))$.  Then we can express $H^0(X,\ce_2)$ 
as the intersection of two maximal $q$-isotropic subspaces and 
use the fact \cite[pp. 735ff]{griffiths-harris} that the 
dimension (mod $2$) of such an intersection is constant under 
deformation. 

To prove that there are no more irreducible components than 
$d_{g,r}(\mathbf{m})$ we first note that when $\mathbf{m}\equiv 
\mathbf{0} \pmod r$, and $g=1$, the claim holds by Lemma 
\ref{elliptic}.  Similarly, if $r=2$ and $g \geq 1$ then we have 
the following lemma. 

\begin{lemma} \label{irred-2}
  If $m_i$ is even for every $i$, and if
  $\overline{\scheme{S}}^{1/2,\bm \text{ even}}_{g,n}$ and 
  $\overline{\scheme{S}}^{1/2,\bm \text{ odd}}_{g,n}$
  are the respective inverse images under $e$ of $0$
  and $1$, then $\overline{\scheme{S}}^{1/2,\bm \text{ even}}_{g,n}$ and
  $\overline{\scheme{S}}^{1/2,\bm \text{ odd}}_{g,n}$ are irreducible.
\end{lemma}

The proof of Lemma~\ref{irred-2} is an easy generalization of 
Cornalba's proof \cite[Lemma 6.3]{cornalba:theta} of the 
corresponding result when $n=0$.  Cornalba's idea is to check the 
result explicitly for the cases $g=1$ and $g=2$ (i.e., write out 
all of the square roots of the canonical bundle---or in our case, 
the square roots of the bundle $\omega(-\sum m_i q_i)$).  And 
then one can use induction and a degeneration argument to reduce 
to the lower genus case.  This argument works just as well for 
our case, except that in the case of $g=1$ and $g=2$, all of the 
$m_i$ must be even to write out the square roots in the form 
needed to see that monodromy acts transitively on the even 
(respectively, odd) square roots. 

The proof of Theorem \ref{disjoint} now follows from the theorem 
below, which is a generalization of  results of Sipe 
\cite{sipe:roots} and Hain \cite[\S13]{hain:mg-torelli}. 

\begin{theorem}\label{transitive}
  For any fixed, smooth curve $X$ of genus $g>0$, if $\ell_{g,r}
  (\mathbf{m})=1$ then the monodromy group acts transitively on the
  set $\sgrnm[X]$ of isomorphism classes of spin structures on $X$
  (where two spin structures that differ by an automorphism of $X$ are
  not considered to be isomorphic).  And if $\ell_{g,r}(\mathbf{m})>1$
  then the monodromy group acts transitively on each fibre of the map
  $\sgrnm [X] \xrightarrow{[r/\ell_{g,r}(\mathbf{m})]}
  \overline{\scheme{S}}^{1/\ell_{g,r},\bm'}_{g,n}[X]$.
\end{theorem}

To prove Theorem~\ref{transitive}, we need to further study the 
 action of the mapping class group on the fibres of $\sgrnm$ over
 $\schememgn$.  In particular, note that if $\Gamma_{g,n}$ is the (pure)
 mapping class group for $n$-pointed curves, and if $\jac_r X$ is the
 group of $r$-torsion points in the Jacobian of a  fixed curve $[X] \in
 \schememgn$, then the group $\Gamma_{g,n}$ acts on the principal homogeneous
 $\jac_r X$-space $\sgrnm[X]$ in a way compatible with the usual
 action of $\Gamma_{g,n}$ on $\jac_r X$.  

The monodromy action induces a homomorphism from $\Gamma_{g,n}$
into the group of affine (with respect to $\jac_r X$), invertible
transformations $\mathcal A$ of $\sgrnm[X]$. And we have the 
following diagram: \begin{diagram} & & & & \Gamma_{g,n} \\ & & & 
& 
\dTo \\ 0 & 
\rTo & \jac_r X & \rTo^{T} & \mathcal A&\rTo & GL_{2g}(\mathbb Z /
r\mathbb Z ) & \rTo & 1\\
\end{diagram} Here the horizontal sequence is exact, and the map $T$
simply takes an element $\eta$ to the automorphism that 
translates by $\eta$.  Of course, since the mapping class group 
preserves the intersection product on $H_1(X,\mathbb Z)$, the 
image of $\Gamma_{g,n}$ in $GL_{2g}(\mathbb Z / r\mathbb Z )$ is 
a subgroup of the symplectic group $SP_{2g}(\mathbb Z / r\mathbb 
Z ).$ And it is well known (see, for example \cite[pg. 
178]{magnus}) that the image of $\Gamma_{g,n}$ is all of 
$SP_{2g}(\mathbb Z / r\mathbb Z ).$ 

We are especially interested in the elements of $\Gamma_{g,n}$ 
which act trivially on the homology $H_1(X,\mathbb Z)$.  The 
subgroup $\mathcal{I}_{g,n}$ of all such elements  is called the 
{\em Torelli group}. The main step in the proof of 
Theorem~\ref{transitive} is the following lemma. 

\begin{lemma}\label{mono-lemma}
  If $g \geq 2$, the image of the Torelli group in the group of translations is the
  subgroup of $\jac_r X$ generated by $\{2, m_1, m_2, \dots ,
  m_n \} \jac_r X$, and if $g=1$ then the image of the Torelli group is generated by 
  $\{m_1, \dots, m_n\} \jac_rX$.  More exactly, if $E$ is the image of $\Gamma_{g,n}$ in
  $\mathcal A$, then we have the following exact sequences.

For $g \geq 2$ $$0 \rTo <2,m_1,\dots,m_n>\jac_r X
\rTo E \rTo SP(2g, \mathbb Z/r \mathbb Z) \rTo 1$$
and for $g=1$ $$0 \rTo <m_1, \dots, m_n>\jac_rX \rTo E \rTo SP(2, 
\mathbb{Z}/r \mathbb{Z}) \rTo 1.$$ 
\end{lemma}

The lemma implies Theorem~\ref{transitive} (transitivity of the
monodromy action) in the case that $\ell_{g,r}(\mathbf{m})=1$.  And on
the other hand, if $\ell_{g,r}(\mathbf{m})>1$, then the map $\sgrnm[X]
\xrightarrow{[r/\ell_{g,r}(\mathbf{m})]} \scheme{S}_{g,n}^{1/\ell,\bm'}
[X]$ is equivariant under the action of the mapping class group, 
and any two spin structures that map to the same point of 
$\scheme{S}^{1/\ell,\bm'}_{g,n}[X] $ must differ by a point of 
$\ell_{g,r}(\mathbf{m}) \cdot \jac_r X$.  Therefore, the orbits 
of the action of $\ell_{g,r}(\mathbf{m}) \cdot \jac_r X$, and 
hence also the orbits of the action of $\Gamma_{g,n}$, are 
exactly the fibres of this map.  Thus all that is necessary for 
the proof of Theorem~\ref{transitive} is to prove the lemma.  The 
following proof generalizes and expands ideas from the results of 
Hain \cite[\S 13]{hain:mg-torelli} and Sipe \cite{sipe:roots}. 

\begin{proof}(of Lemma \ref{mono-lemma}) Rather than studying the $r$th roots of the bundle
  $\omega_X (-\sum m_i p_i)$, it is convenient to dualize and study
  the action on the $r$th roots of the bundle $\cl := TX (\sum m_i
  p_i)$, where $TX$ is the tangent bundle to $X$, Let $\cl^{o}$ denote
  the $S^1$-bundle, obtained by removing the zero section from $\cl$
  and retracting to the unit circle in each fibre.  The Euler class of
  $\cl^{o}$ is $2-2g+\sum m_i$ and is divisible by $r$, so the end of
  the Gysin sequence gives a short exact sequence
\begin{equation}0 \rightarrow 
\mathbb Z/r \mathbb Z \rightarrow H_1(\cl^{o}, \mathbb Z/ r \mathbb Z)
\rightarrow H_1(X , \mathbb Z/ r \mathbb Z) \rightarrow
0. \label{gysin} \end{equation} If $\sigma: [0,1] \rightarrow X$ is a
loop in $X$ that does not pass through any of the $n$ points $p_i$,
then there is a canonical lift of $\sigma$ to a loop $\tilde \sigma$
in $\cl^{o}$, which is defined as $\displaystyle \tilde \sigma (t) =
\left( \sigma (t) , \frac{\dot\sigma (t)}{|| \dot \sigma(t)||}\right)
$.  Note that for any small loop $\lambda$, homotopic to zero, which
does not pass through any of the $p_i$, the lift $\tilde \lambda$
corresponds to $1$ or $-1$ in the group $\mathbb Z/r \mathbb Z$ on the
left side of the Gysin sequence~\eqref{gysin}.

Given any $r$th root $\cn$ of $\cl$, there is a natural covering 
map of $S^1$-bundles $p: \cn^{o} \rightarrow \cl^{o}$, which 
induces a map $\pi_1(\cl^{o}) \rightarrow \pi_1(\cl^{o})/p_* 
\pi_1(\cn^{o}) \cong 
\mathbb Z/ r \mathbb Z $, and thus a homomorphism $H_1 (\cl^{o})
\rightarrow \mathbb Z/ r \mathbb Z$.  Moreover, this map takes a lift
$\tilde \lambda$ of a small loop $\lambda$, homotopic to zero in 
$X$, and maps it to $\pm 1$.  In other words, the $r$th root 
induces a splitting of the Gysin sequence.  Conversely, given any 
such splitting, there is a covering space of $\cl^{o}$ that can 
easily be seen to define an $r$th root of $\cl$.  In particular, 
we have the following fact, originally due to Sipe 
\cite{sipe:roots} in the case when $n=0$.  \label{root-split} 
\begin{fact}
 The $r$th roots of $\cl^{o}$ are in
a natural one-to-one correspondence with the splittings of the
sequence \eqref{gysin}.
\end{fact}

We will compute the action of certain elements of $\Gamma_{g,n}$ 
on $H_1(\cl^{o},\mathbb Z/r \mathbb Z)$.  Choose a basis for 
$H_1(\cl^{o},\mathbb Z/r \mathbb Z)$ in the following way.  First 
take an explicit choice of cycles $B_1, B_2, \dots, B_{2g}$ which 
form a basis of $H_1(X,\mathbb Z/r \mathbb Z)$, and such that for 
all $i,j,k$, the intersection form $<B_i , B_{i+1}> = 1$, and 
$<B_j,B_k>=0$ if $k$ is not $j+1$ or $j-1$.  An example of such a 
collection of cycles is depicted in Figure~\ref{basis}. 

\begin{figure}
\includegraphics{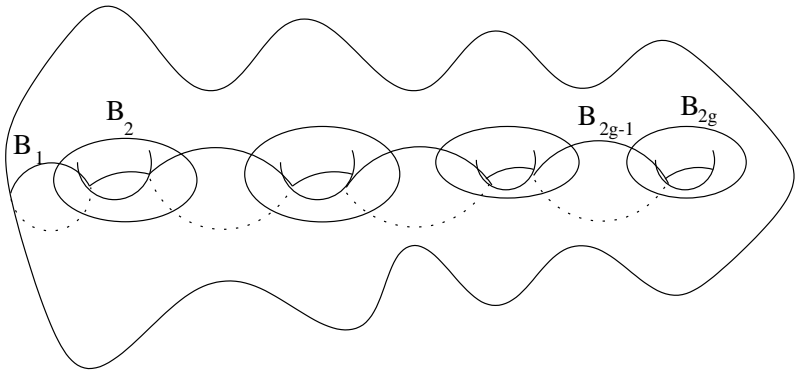}
\caption{\label{basis}
A collection $\{ B_1, \dots , B_{2g} \}$ of cycles forming a basis for 
$H_1(X,\mathbb Z)$, such that the
intersection form \mbox{$<B_i , B_{i+1}> = 1$}, and the intersection
numbers for all other pairs are zero.}
\end{figure}

We also insist that none of these cycles contain any of the $n$ marked
points of $X$.  Now lift these basis cycles in the usual way to
$\tilde B_1, \dots \tilde B_{2g}$.  As before, take a small loop
$\lambda$, which does not pass through any of the $p_i$ and is
homotopic to zero, and lift it to $\tilde \lambda$.  Now the set
$\{\tilde B_1, \dots, \tilde B_{2g}, \tilde \lambda\}$ gives a basis
for $H_1(\cl^{o},\mathbb Z/r \mathbb Z)$, and it defines a splitting
of the sequence \eqref{gysin}.  Any other splitting is given by simply
adding this splitting map $H_1(X) \rightarrow H_1(\cl^{o})$ to an
element of $H^1(X,\mathbb Z/ r \mathbb Z) \cong \jac_r X$. And given a
splitting induced by $\{\tilde B_1, \dots, \tilde B_{2g}, \tilde
\lambda\}$, we will index the remaining $r$th roots by the
corresponding $1$-cocycle.

Given any simple closed curve $\gamma$ on $X$ that misses the 
marked points $\{p_i\}_{i=1}^n$, there is a corresponding Dehn 
twist $T_{\gamma}$ that is an element of $\Gamma_{g,n}$.  Using 
the Picard-Lefschetz Theorem \cite{zariski:alg-surf}, one can 
show (see, for example \cite[Proposition 3.1]{sipe:roots}) that 
the action of $T_{\gamma}$ on the lift $\tilde \sigma$ to 
$\cl^{o}$ of a cycle $\sigma$ on $X$ is given by the following 
generalized Picard-Lefschetz formula. 
\begin{equation} \label{PL} 
T_{\gamma}(\tilde \sigma) = \tilde \sigma + <\sigma , \gamma>
\tilde \gamma \end{equation} In particular, the action of $T_{\gamma}$
 on $\tilde \lambda$ is trivial.

\begin{figure}
\includegraphics{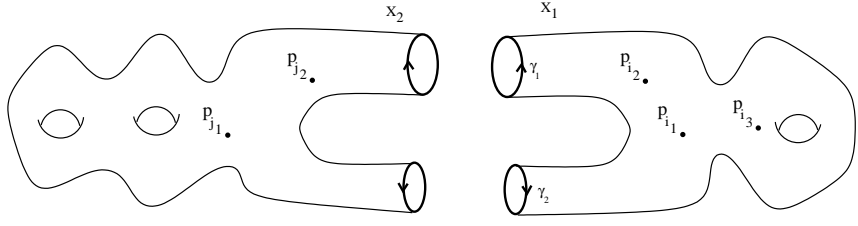}
\caption{\label{bpscc}
A bounding pair of simple closed curves, separating $X$ into
$X_1$ and $X_2$ of genus $g_1=1$ and $g_2=2$ respectively, and containing
the marked points $p_{i_1}, \dots, p_{i_m}$ and the points
$p_{j_1},\dots,p_{j_{n-m}}$, respectively.}
\end{figure}

Now consider a cycle consisting of a bounding pair (sometimes 
called a \emph{cut pair}) of simple closed curves $\gamma_1$ and  
$\gamma_2$ in $X$ which separate $X$ into two surfaces $X_1$ and 
$X_2$ and such that the cycle $\gamma_1 + \gamma_2$ is homologous 
to zero.  Assume these two surfaces are of genus $g_1$ and $g_2$, 
respectively, and assume that the points $p_{i_1}, 
\dots, p_{i_m}$ are in $X_1$, and the remaining points are in $X_2$ as
in Figure~\ref{bpscc}. If the cycles $\gamma_1$ and $\gamma_2$ have
non-trivial intersection number with $B_i$, then $<B_i , \gamma_1> =
<B_i , \gamma_2>$, and the generalized Picard-Lefschetz formula
(\ref{PL}) gives
\begin{equation}
T_{\gamma_1} \circ T_{\gamma_2} (\tilde B_i) =T_{\gamma_2} \circ
T_{\gamma_1} (\tilde B_i) = \tilde B_i + <B_i , \gamma_1>(\tilde
\gamma_1 + \tilde \gamma_2).
\end{equation}
Since the cycle $\gamma_1+\gamma_2$ is homologous to zero, the  
cycle $\tilde 
\gamma_1 +\tilde \gamma_2$ in $H_1(\cl^{o},\mathbb Z/r \mathbb Z)$  is an integral 
multiple---say $a_{\gamma_1 + 
\gamma_2}$---of $\tilde \lambda$.  Thus the action induced by 
$T_{\gamma_1} \circ T_{\gamma_2}$ on the set of $r$th roots of 
$\cl$ is simply the one that translates $H^1(X,\mathbb Z/ r 
\mathbb Z)$ by the cocycle that takes $B_i$ to $a_{\gamma_1 + 
\gamma_2} \cdot <B_i , \gamma_1>$.

We need, therefore, to compute $a_{\gamma_1 + \gamma_2}$.  To 
this end, define a $2$-chain in $\cl^{o}$, using, as in 
\cite[Proposition 3.2]{sipe:roots}, a singular unit vector field 
on $X_1$.  To obtain this vector field, simply glue together 
$2g_1$ copies of the field represented in Figure~\ref{vectfield} 
to make a vector field on the $X_1$. 

\begin{figure}
\includegraphics{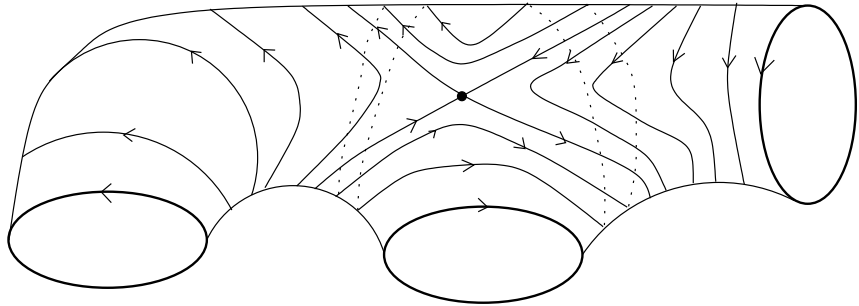}
\caption{\label{vectfield}
A singular vector field, which when joined to $2g_1$ other such
pieces, defines a singular vector field in $TX^o$, homeomorphic 
to $X_1$.  Removing a small disc near each singularity gives a 
$2$-chain in $\cl^{o}$ with boundary $\displaystyle \tilde 
\gamma_1 + \tilde \gamma_2$ plus $2g_1 - \sum_{p_{i_k} \in X_1} 
m_{i_k}$ $1$-cycles homologous to $\tilde \lambda$.} 
\end{figure}

Cutting out a small disc around each point where the vector field is
singular, we obtain a surface in the unit tangent bundle $TX^{o}$
bounded by $\tilde \gamma_1$, $\tilde \gamma_2$, and $2g_1$ cycles
homologous to $\tilde \lambda$.  If we consider the vector field as a
section of $\cl^{o}$ instead, singularities at the marked points
$p_{i_k}$ do not require the removal of the corresponding discs,  and
the $2$-chain in $\cl^{o}$ is bounded by $\tilde \gamma_1$, $\tilde
\gamma_2$, and by $\displaystyle 2g_1 - \sum_{p_{i_k} \in X_1}
m_{i_k}$ cycles homologous to $\tilde \lambda$.  Thus the action
induced by $T_{\gamma_1} \circ T_{\gamma_2}$ on the set of $r$th
roots of $\cl$ is translation by the cocycle that takes $B_i$ to 
$ 
\left ( \displaystyle 2g_1 - \sum_{p_{i_k} \in X_1} m_{i_k}\right )
\cdot <B_i , \gamma_1>$, i.e., it is $\displaystyle\left ( 2g_1 -
\sum_{p_{i_k} \in X_1} m_{i_k}\right )$ times the Poincar\'e dual of
$\gamma_1$.  Note that this is well-defined because the sum
$\displaystyle 2g_1 +2g_2 - \sum_{i} m_{i} = 2g-2 - \sum m_i$ is
divisible by $r$.

It is easy to see that for any basis cycle $B_i$, for any 
integers $g_1$ and $g_2$ that sum to $g-1$, and for any subset 
$\{ p_{i_k} \}$ of the marked points, there is a bounding pair 
$\gamma_1$ and $\gamma_2$ that are each dual to $B_i$, and which 
cut $X$ into two surfaces $X_1$ and $X_2$ of genera $g_1$ and 
$g_2$ respectively. Moreover, this can be done so that $X_1$ 
contains the points $\{ p_{i_k} \}$, and $X_2$ contains the 
remaining marked points. Thus the action of $\Gamma_{g,n}$ on the 
$r$th roots of $\cl$ includes all translations by elements in the 
subgroup $<2, m_1, \dots, m_n> H^1 (X,\mathbb{Z}/r \mathbb{Z})$ 
if $g>1$ and $<m_1, \dots, m_n>H^1 (X, \mathbb{Z}/r\mathbb{Z})$ 
if $g=1$. 

 Proposition~\ref{Torelli} below shows that the Torelli group is
 generated
 by Dehn twists on bounding 
 (cut) pairs, as above, and by Dehn twists on separating simple closed
 curves (also called bridges).  Thus  all that remains in the proof of the
 lemma is to show that Dehn twists along separating simple closed
 curves act trivially on $\sgrnm[X]$.  But this follows immediately
 from the generalized Picard-Lefschetz formula~\eqref{PL} and the fact
 that the intersection of a separating simple closed curve with the
 $B_i$ is zero.  This completes the proof of the lemma and
 also the proof of Theorem~\ref{transitive}.
\end{proof}

The following proposition gives a set of generators for the 
Torelli group in all genera greater than zero, and for all $n$.  
Although this result follows easily from known results, it does 
not seem to exist, as such, in the literature, so we include it 
here for completeness. 

\begin{proposition}\label{Torelli}
For $g \geq 1$ and $n \geq 0$ the Torelli group 
$\mathcal{I}_{g,n}$ is generated by Dehn twists on bounding pairs 
of simple closed curves and Dehn twists on separating simple 
closed curves.  If $g 
\neq 2$, a finite subset of such twists generate $\mathcal{I}_{g,n}$
\end{proposition} 

\begin{proof}
In the case that $g \geq 3$ and $n=0$ the result is due to Powell 
\cite{powell:torelli}, and Johnson \cite{johnson:torelli-1} 
proved that a finite subset of these will suffice to generate 
$\mathcal{I}_{g,0}$.  In the case of $g=1$ and $n=0$ or $n=1$ it 
is a classical fact that $\mathcal{I}_{1,0} 
=\mathcal{I}_{1,1}=\{1\}$.  And for $g=2$ and 
$n=0$ the result of the proposition was proved by G. Mess 
\cite{mess:torelli}, but no finite set of generators is known. 

For every $n\geq 1$ Birman \cite[pp158-60]{birman:blmcg} gives a 
finite collection of bounding pairs of simple closed curves 
(unlike those of Johnson, these bounding pairs may cut the 
surface so that one of the resulting pieces has genus zero) such 
that Dehn twists on these generate the kernel $K_{g,n}$ of the 
natural homomorphism $i_*:\Gamma_{g,n} \rTo \Gamma_{g,n-1}$, and 
since $i_*$ is surjective, we have an exact sequence $$1 \rTo 
K_{g,n} 
\rTo \mathcal{I}_{g,n} \rTo \mathcal{I}_{g,n-1} \rTo 1.$$  By  
induction, the Proposition holds for all $n$. 
\end{proof}

\begin{nb}
  For $g \geq 2$ the number of elements in
  ${\mbox{$\sgrnm$}}^{even}[X]$ and ${\mbox{$\sgrnm$}}^{odd}[X]$ is
  easily computed as the number of even, respectively odd,
  theta-characteristics times the degree of the morphism
  $[r/\ell_{g,r} (\mathbf{m})]$ (which is the order of the group
  $2\cdot \jac_r X$).  And it is known \cite{cornalba:theta} that of
  the $2^{2g}$ square roots of $\omega_X$, $2^{g-1}(2^{g} + 1)$ are
  even, and $2^{g-1}(2^{g} - 1)$ are odd; so there are $ r^{2g}(1/2 +
  1/{2^{g+1}})$ even elements in $\scheme{S}^{1/r,\bm}_{g,n}[X]$ and
  $r^{2g}(1/2 - 1/{2^{g+1}})$ odd, (Some of these elements may be
  identified by automorphisms of $X$).
  
  Similarly, for $g=1$ the number of elements in the component of
  $\scheme{S}^{1/r,\mathbf{0}}_{1,n}$ indexed by $d$ is well-known
  \cite[11\S 3 Prop 3.6 (2)]{husemoller:elliptic-curves}, and easily
  calculated to be $(r/d)^2 \prod_{p|(r/d)} (1-1/p^2)$ when $X$ has no
  automorphisms (generically when $n\geq2$).  And the degree of
  $[r/\ell_{1,r}(\mathbf{m})]$ is $(r/\ell_{1,r}(\mathbf{m}))^2$.
\end{nb}

\begin{nb} Using the techniques of the  proof of
  Lemma~\ref{mono-lemma} to calculate the action of various Dehn
  twists on $\sgrnm[X]$, it is not hard to see that, when $r>2$, the
  stabilizer of one point in $\sgrnm[X]$ is not normal in
  $\Gamma_{g,n}$.  That is, the finite cover $\sgrnmbar \rightarrow
  \schememgnbar$ is not regular.  This is an interesting contrast to
  the fact that the subgroup of $\Gamma_{g,n}$ which fixes {\em all}
  the points of $\sgrnm[X]$ is normal \cite[Theorem C]{sipe:roots}.
\end{nb}

\subsubsection{Irreducibility over General Ground Fields}

  Although in the previous section we have been working  over the
  complex numbers $\mathbb C$, $r$-spin curves form a smooth 
  Deligne-Mumford stack over $\mathbb Z[1/r]$.  Moreover, as in the case of
  $\mgbar$ (see \cite[\S5]{deligne-mumford}), irreducibility in
  characteristic $0$ gives the same result in any characteristic
  relatively prime to $r$.

\begin{theorem}\label{all-char}
  For any field $k$ with characteristic prime to $r$, the moduli space
  $\sgrnm$ and and its compactification $\sgrnmbar$ of pure spin
  curves over $k$ are irreducible if $\ell_{g,r}(\mathbf{m})=1$.  
  In general, over any algebraically closed field, the moduli space is
  the disjoint union of $d_{g,r}(\mathbf{m})$ irreducible components.
\end{theorem}

\begin{proof}
The proof is essentially the same as for the case of $\schememg$, 
but we will sketch the main steps. 

%fix notation
First, note that it is well-known and straightforward to prove 
\cite[Lemma~2.3]{vistoli} that an algebraic stack is irreducible 
or connected if and only if its coarse moduli space has the same 
property.   It is easy to see from the universal deformation of 
stable spin curves that the stack $\stacksgrnm$ of smooth spin 
curves is an open dense substack of the stack $\stacksgrnmbar$ of 
stable spin curves, and thus its irreducible components are the 
non-empty intersections of $\stacksgrnm$ with the irreducible 
components of $\stacksgrnmbar$ \cite[Prop 4.15]{deligne-mumford}.  
So it suffices to prove the theorem in the case of 
$\stacksgrnmbar$. 

Since $\stacksgrnmbar$ is smooth, its  connected components are 
irreducible \cite[Prop 4.16]{deligne-mumford}.  Moreover, if for 
each $s$ in $\spec \mathbb Z[1/r]$, we define $n(s)$ to be the 
number of connected components of the geometric fibre of 
$\stacksgrnmbar$ over $s$, then $n(s)$ is a constant function 
\cite[Prop 4.17]{deligne-mumford}. 

And finally, by Theorem~\ref{transitive} the assertion of
Theorem~\ref{all-char} holds over $\mathbb C$ for the moduli 
space $\sgrnmbar$, and thus for 
 the stack 
$\stacksgrnmbar$ .  Consequently the assertion must hold for both 
the stack and its moduli space over all algebraically closed 
fields. And in the case that $\ell_{g,r}(\mathbf{m})=1$ the stack 
and its moduli space must be irreducible over any field. 
\end{proof}

\section*{Acknowledgments}

I wish to thank Dan Abramovich, Steve Humphries, and Arkady 
Vaintrob for helpful discussions and suggestions regarding this 
work.  I am also grateful to Heidi Jarvis for help with 
typesetting. 

%\nocite{cgm:r-torsion}
%\bibliographystyle{amsplain}
%\bibliography{merged}

\providecommand{\bysame}{\leavevmode\hbox to3em{\hrulefill}\thinspace}

\end{document}